\documentclass[12pt,twoside]{article}

\usepackage[T1]{fontenc}
\usepackage[latin1]{inputenc}
\usepackage[english]{babel}
\usepackage[babel]{csquotes}

\usepackage{cite}

\usepackage{amssymb}
\usepackage{amsmath}
\usepackage{amsthm}
\usepackage{latexsym}
\usepackage{graphicx}
\usepackage{mathrsfs}
\usepackage{mathrsfs}
\usepackage{bbm}

\usepackage[usenames,dvipsnames]{color}

\theoremstyle{plain}
\newtheorem{thm}{Theorem}[section]
\newtheorem{theorem}{Theorem}[section]

\newtheorem{proposition}[thm]{Proposition}

\newtheorem{definition}[thm]{Definition}
\theoremstyle{definition}

\newtheorem{remark}[thm]{Remark}

\def\enne{\mathbb{N}}
\def\zeta{\mathbb{Z}}

\def\erre{\mathbb{R}}

\def\P{\mathbb{P}}

\def\E{\mathop{{}\mathbb{E}}}
\def\cL{\mathscr{L}}
\def\cF{\mathscr{F}}
\def\cB{\mathscr{B}}
\def\eps{\varepsilon}

\def\beq{\begin{equation}}
\def\eeq{\end{equation}}

\def\to{\rightarrow}
\def\wto{\rightharpoonup}
\def\wstarto{\stackrel{*}{\rightharpoonup}}
\def\embed{\hookrightarrow}
\def\cembed{\stackrel{c}{\hookrightarrow}}

\def\norm #1{\left\|#1\right\|}

\def\sp #1#2{\left<#1,#2\right>}
\newcommand\ip\sp

\title{\huge\rm Optimal distributed control\\of a stochastic Cahn-Hilliard equation
			\footnote{
			This paper was 
			funded by Vienna Science and Technology Fund (WWTF) 
			through Project MA14-009.}
		\\[1cm]}

\author{	{\large\sc Luca Scarpa}\\
		{\normalsize E-mail: \texttt{luca.scarpa@univie.ac.at}}\\[.5cm]
		{\small Faculty of Mathematics, University of Vienna}\\
		{\small Oskar-Morgenstern-Platz 1, 1090 Vienna, Austria}
	}

\date{}

\begin{document}

\maketitle  

\begin{abstract}
  We study an optimal distributed control problem associated to a stochastic 
  Cahn-Hilliard equation with a classical double-well potential and Wiener 
  multiplicative noise, 
  where the control is represented by a source-term in the definition of the chemical potential.
  By means of probabilistic and analytical compactness arguments, 
  existence of an optimal control is proved. Then the linearized system and the 
  corresponding backward adjoint system are analysed through 
  monotonicity and compactness arguments,
  and first-order necessary conditions for optimality are proved.
  \\[.5cm]
  {\bf AMS Subject Classification:} 35K55, 35R60, 60H15, 80A22, 82C26.\\[.5cm]
  {\bf Key words and phrases:} stochastic Cahn-Hilliard equation, phase separation, optimal control,
  linearized state system, adjoint state system, 
  first-order optimality conditions.
\end{abstract}

\pagestyle{myheadings}
\newcommand\testopari{\sc Luca Scarpa}
\newcommand\testodispari{\sc Optimal control of a stochastic CH equation}
\markboth{\testodispari}{\testopari}


\thispagestyle{empty}

\section{Introduction}
\setcounter{equation}{0}
\label{sec:intro}

The pure Cahn-Hilliard equation on a smooth bounded domain $D\subset\erre^N$, $N=2,3$,
can be written in its simplest form as
\[
  \partial_t y - \Delta w = 0\,, \quad w=-\Delta y + \Psi'(y) - u \qquad\text{in } (0,T)\times D\,,
\]
where $T>0$ is a fixed final time, $y$ and $w$ denote the order parameter and the 
chemical potential of the system, respectively, and $u$ represents a given distributed source term.
Furthermore, $\Psi'$ is the derivative of a so-called double-well potential $\Psi$, which may be
seen as the sum of a convex function and a concave quadratic perturbation:
typical examples of $\Psi$ which are relevant in applications are discussed in \cite{col-gil-spr}.
Usually, in order to ensure the conservation of the mean on $D$,
the equation is complemented by homogenous Neumann conditions 
for both $y$ and $w$, and a given initial value, namely
\[
  \partial_{\bf n} y = \partial_{\bf n} w = 0 \quad\text{in } (0,T)\times\partial D\,,
  \qquad y(0)=y_0 \quad\text{in } D\,,
\]
where $\bf n$ denotes the outward normal unit vector on $\partial D$.

The Cahn-Hilliard equation was originally introduced in \cite{cahn-hill}
(see also \cite{ell-st, ell-zhen, novick-cohen}) to capture the spinodal decomposition
phenomenon occurring in a phase-separation of a binary metallic alloy.
The mathematical literature on the deterministic Cahn-Hilliard equation
has been widely developed in the last years, especially 
in much more general settings as the presence of
viscosity terms and dynamic boundary conditions:
in this direction we mention, among all, the contributions
(as well as the references therein)
\cite{bcst1, col-fuk-eqCH, col-gil-spr, cher-gat-mir, cher-mir-zel, cher-pet, colli-fuk-CHmass,
gil-mir-sch, mir-sch, scar-VCHDBC} on the well-posedness of the system, and
\cite{col-fuk-diffusion, col-scar, gil-mir-sch-longtime}
on asymptotic behaviour of the solutions.
Optimal distributed and boundary control problems have been studied
in the context of Allen-Cahn and Cahn-Hilliard equations
in the works \cite{col-far-hass-gil-spr,
col-far-hass-spr, col-gil-pod-spr, col-gil-spr-contr3,
col-gil-spr-contr, col-gil-spr-contr2, hinter-weg, roc-spr}.

More recently, in order to account also for the random vibrational 
movements at a microscopic level in the system, which may be 
of magnetic, electronic or configurational nature, 
the equation has been modified by adding a cylindrical Wiener process $W$
(see \cite{cook, lee-CH}). This has resulted in the well-accepted version
of the stochastic Cahn-Hilliard equation
\begin{align}
  \label{eq:1}
  dy - \Delta w\,dt = B(y)\,dW \qquad&\text{in } (0,T)\times D=:Q\,,\\
  \label{eq:2}
  w=-\Delta y + \Psi'(y) - u \qquad&\text{in } (0,T)\times D\,,\\
  \label{eq:bound}
  \partial_{\bf n} y = \partial_{\bf n} w = 0 \qquad&\text{in } (0,T)\times\partial D=:\Sigma\,,\\
  \label{eq:init}
   y(0)=y_0 \qquad&\text{in } D\,,
\end{align}
where $B$ is a stochastically integrable operator with respect to $W$.
The mathematical literature on the stochastic Cahn-Hilliard and 
Allen-Cahn equations is significantly less developed. Let us mention the works
\cite{corn, daprato-deb, elez-mike, scar-SCH}
dealing with existence, uniqueness and regularity
for the pure equation, and \cite{hao-stochVCH, ju-stochVCH, scar-SCH}
for an analysis of the viscous case
in terms of well-posedness, regularity and vanishing viscosity limit.
We point out for completeness also the contributions
\cite{ant-kar-mill} for a study of a stochastic Cahn-Hilliard 
equation with unbounded noise, and \cite{deb-zamb, deb-goud, goud} dealing
with stochastic Cahn-Hilliard equations with reflections.
The reader can refer also to \cite{bauz-bon-leb, orr-scar}
for the context of stochastic Allen-Cahn equations, and
\cite{feir-petc} for a study of a diffuse interface model 
with termal fluctuations.

While the literature on stochastic optimal control
problems is widely developed, we are not aware of any result dealing 
with controllability of the stochastic Cahn-Hilliard equation. 
The main novelty of the present contribution is to provide
a first study in this direction, and represents
a starting point for the study of optimal control problems
associated to the wide class of more general phase-field models
with stochastic perturbation.
Optimal control problems have been studied in the stochastic case 
especially in connection with the stochastic maximum principle:
the reader can refer to \cite{yong-zhou} for a general treatment
on the subject. 
Let us mention the works \cite{orr-contr, orr-fur} dealing with 
stochastic maximal principle for nonlinear SPDEs with 
dissipative drift, \cite{du-meng} for optimal control
of stochastic evolution equations in Hilbert spaces, and
\cite{fur-hu-tess} for
a study on optimal control of SPDEs with 
control contained both in the drift and the diffusion.
Let us also point out \cite{orr-vev} for a stochastic optimal control
problem on infinite time horizon, \cite{orr-tess-vev}
on ergodic maximum principle,
and
\cite{brzez-serr} for a study on optimal 
relaxed controls of dissipative SPDEs. Stochastic optimal 
control problems have also been considered in \cite{barbu-rock-contr}
in the context of the Schr\"odinger equation.

In the present contribution we are interested in studying a distributed optimal 
control problem associated to the stochastic pure Cahn-Hilliard equation,
where the control is the source term $u$ in the definition of the chemical potential
and the cost functional is of standard quadratic tracking-type. More precisely, 
we want to minimize
\beq\label{cost}
  J(y,u):=\frac{\alpha_1}{2}\E\int_Q|y-x_Q|^2 + 
  \frac{\alpha_2}{2}\E\int_D|y(T)-x_T|^2 +
  \frac{\alpha_3}{2}\E\int_Q|u|^2\,,
\eeq
subject to the state equation \eqref{eq:1}--\eqref{eq:init} and
a control constraint on $u\in \mathcal U$, where
$\mathcal U$ is a suitable convex closed subset of
$L^2(\Omega\times Q)$
which will be specified in Section~\ref{sec:main} below.
Here, $\alpha_1$, $\alpha_2$, $\alpha_3$ are nonnegative constants,
$x_Q$ and $x_T$ are given functions in $L^2(\Omega\times Q)$
and $L^2(\Omega\times D)$, respectively.
The main results of this work are the existence of a relaxed optimal control
and the proof of first-order necessary conditions for optimality.

The first step of our analysis consists in studying the control-to-state mapping.
In particular, we show that for every admissible control $u\in\mathcal U$, 
the state system \eqref{eq:1}--\eqref{eq:init} admits a unique solution $y$,
and the map $S:u\mapsto y$ is Lipschitz-continuous in some suitable spaces.
Consequently, the cost functional $J$ can be expressed in a reduced form
only in terms of the control $u$, i.e.~introducing the reduced cost functional $\tilde J$ as
\[
\tilde J(u):=J(S(u), u)\,, \qquad u\in\mathcal U\,.
\]

At this point, in the deterministic setting
the most natural necessary condition for optimality of $\bar u\in \mathcal U$ would read
\[
  D\tilde J(\bar u)(v-\bar u)\geq 0 \qquad\forall\, v\in\mathcal U\,,
\]
where $D\tilde J$ represents the derivative of $\tilde J$
at least in the sense of G\^ateaux.
In this direction,
the classical approach consists in showing that 
the map $S$ is Fr\'echet-differentiable, hence so is 
$\tilde S$ by the usual chain rule for Fr\'echet-differentiable functions,
and to characterize the derivative $DS(\bar u)$ as 
the solution of a suitable linearized system.
In the context of Cahn-Hilliard equations with possibly degenerate potentials
(for example if $\Psi$ is the double-well logarithmic potential), the
Fr\'echet differentiability of the control-to-state mapping is usually obtained
by requiring sufficient conditions in the box constraint for $u$,
ensuring at least that $\Psi''(\bar y) \in L^\infty(Q)$, where $\bar y:=S(\bar u)$
and $\Psi''$ is the second derivative of $\Psi$
(for example that $\mathcal U$ is contained in a closed ball in $L^\infty(Q)$).

However, if we add a stochastic perturbation in the equation, 
under reasonable assumptions on the data it is {\em not} possible to prove
that $\Psi''(\bar y)$ is uniformly bounded in $L^\infty(\Omega\times Q)$, even 
if we add a constraint on the $L^\infty$-norm in the definition of the admissible controls.
This behaviour gives rise to several nontrivial difficulties: 
among all, it is not true {\em a priori} that the control-to-state map $S$
is Fr\'echet-differentiable in some space.
This issue is usually overcome in the stochastic setting 
using specific time-variations on the control (the so-called
``spike-variation'' technique). In our case, however, 
we are able to avoid such procedure by analysing explicitly 
the linearized system.
More specifically, we prove that the linearized system admits
a unique variational solution by means of compactness and monotonicity arguments.
Then, we show that the control-to-state mapping is 
G\^ateaux differentiable in a suitable weak sense, and 
that the (weak) G\^ateaux derivative of $S$ can still be identified 
as the unique solution $z$ to the linearized system.
Performing usual first-order variations around a fixed optimal control $\bar u$,
we then prove that the weak G\^ateaux-differentiability is enough to 
ensure first-order necessary conditions for optimality.

The second main issue that we tackle in this work consist in
removing the dependence on $z$ in the first-order necessary conditions
by studying the adjoint problem. As it is well-known,
in the stochastic framework the adjoint problem becomes 
a backward stochastic partial differential equation (BSPDE) of the form
\begin{align}
  \label{eq:1_ad}
  \tilde p=-\Delta p \qquad&\text{in } Q\,,\\
  \label{eq:2_ad}
  -dp -\Delta\tilde p\,dt + \Psi''(y)\tilde p\,dt
  =\alpha_1(y-x_Q)\,dt
  +DB(y)^*q\,dt -q\,dW \qquad&\text{in } Q\,,\\
  \label{eq:bound_ad}
  \partial_{\bf n} p = \partial_{\bf n} \tilde p = 0 \qquad&\text{in } \Sigma\,,\\
  \label{eq:fin_ad}
   p(T)=\alpha_2(\bar y(T)-x_T) \qquad&\text{in } D\,,
\end{align}
where the unknown is the triple $(p,\tilde p, q)$.
Since $\Psi''(y)$ does {\em not} belong to $L^\infty(\Omega\times Q)$,
as we have pointed out above, the adjoint problem cannot be 
framed in any available existence theory for BSPDEs, and
is absolutely nontrivial and interesting on its own.
Through a suitable approximation involving a truncation on $\Psi''$ and a passage to the limit, 
we show existence and uniqueness of a solution to the adjoint problem. 
Furthermore, we prove a suitable duality relation between $z$ and $\tilde p$,
which allows us to express the first-order optimality conditions only in terms
of $\tilde p$ and $\bar u$ in a much more simplified form.

The paper is organized as follows. In Section~\ref{sec:main}
we fix the assumptions, the general setting of the work and the main results.
Section~\ref{sec:wp} contains the proof of well-posedness 
of the state system.
In Section~\ref{sec:exist} we prove that a relaxed optimal control always exists,
using Prokhorov and Skorokhod theorems and natural lower semicontinuity results.
In Section~\ref{sec:map} we study the control-to-state map:
we show that it is well-defined and differentiable in a certain weak sense,
and we identify its (weak) derivative as the unique solution to the linearized problem.
Finally, in Section~\ref{sec:ad} we study the adjoint problem and 
prove the first-order necessary conditions for optimality.

\section{Main results}
\label{sec:main}
\setcounter{equation}{0}

Throughout the paper $(\Omega,\cF, (\cF_{t})_{t\in[0,T]}, \P)$ denotes a filtered 
probability space satisfying the usual conditions, with $T>0$ fixed,
and $W$ is a cylindrical Wiener process on a separable Hilbert space $U$.
The progressive $\sigma$-algebra on $\Omega\times[0,T]$ is denoted by $\mathcal P$.
Furthermore, $D\subset\erre^N$, with $N=2,3$, is a smooth bounded connected domain,
and we use the notation $Q:=(0,T)\times D$ and $Q_t:=(0,t)\times D$ for every $t\in(0,T)$.

For every Hilbert spaces $E_1$ and $E_2$ we denote by $\cL(E_1, E_2)$
and $\cL^2(E_1,E_2)$ the spaces of linear continuous and Hilbert-Schmidt 
operators from $E_1$ to $E_2$, respectively.
The symbols for norms and dualities are $\norm{\cdot}$ and $\ip{\cdot}{\cdot}$,
respectively, with a sub-script indicating the specific spaces in consideration.
We shall use the symbols $\to$, $\wto$ and $\wstarto$
to denote strong, weak, and weak* convergences, respectively.
For any Banach space $E$ and $p\in[1,+\infty]$ 
we shall use the symbols $L^p(\Omega; E)$ and $L^p(0,T; E)$ for the usual 
spaces of Bochner-integrable functions, and the symbols
$C^0([0,T]; E)$ and $C^0_w([0,T]; E)$ for the spaces of
continuous functions
from $[0,T]$ to $E$ endowed with the norm topology or
weak tolopogy, respectively.
If $p,q\in[1,+\infty)$
we shall denote by $L^p_{\mathcal P}(\Omega;L^q(0,T; E))$ the
space of $E$-valued progressively measurable processes $X$ such that
$\E\left(\int_0^T\norm{X(s)}_E^q\,ds\right)^{p/q}<+\infty$.

We define the functional spaces
\[
  H:=L^2(D)\,, \qquad V:=H^1(D)\,, \qquad 
  Z:=\left\{\varphi\in H^2(D):\;\partial_{\bf n}\varphi=0 \text{ a.e.~on } \partial D\right\}\,,
\]
endowed with their natural norms. In the sequel $H$ is identified to $H^*$,
so that $(V,H,V^*)$ is a Hilbert triplet, with dense, continuous and compact inclusions.
The Laplace operator with Neumann homogeneous conditions will be intended 
in the usual variational way as the operator
\[
  -\Delta: V\to V^*\,, \qquad
  \ip{-\Delta x}{\varphi}_V:=\int_D\nabla x\cdot\nabla\varphi\,, \quad x,\varphi\in V\,,
\]
or
\[
  -\Delta: H\to Z^*\,, \qquad
  \ip{-\Delta x}{\varphi}_Z:=-\int_Dx\Delta\varphi\,, \quad x\in H\,,\;\varphi\in Z\,.
\]
We recall also that in the context of Cahn-Hilliard equations it is useful to introduce
the operator $\mathcal N$ as the inverse of $-\Delta$ restricted to the 
subspace of null-mean elements in $V$. More specifically, 
if we denote $x_D:=\frac1{|D|}\ip{x}{1}_V$ for any $x\in V^*$, 
by the Poincar\'e inequality we know that 
\[
  -\Delta: \left\{x\in V: x_D=0\right\} \to \left\{x\in V^*: x_D=0\right\}
\]
is an isomorphism, hence its inverse $\mathcal N$ is well-defined.
Furthermore, it is well-known (see \cite[pp.~979-980]{col-gil-spr}) that
\[
x\mapsto \norm x_*:=\norm{\nabla\mathcal N(x-x_D)}_H + |y_D|\,, \quad x\in V^*\,,
\]
defines a norm on $V^*$, equivalent to the usual one, such that
\beq\label{comp_ineq}
  \forall\sigma>0\,,\quad\exists\,C_\sigma>0: \quad
  \norm x_H^2 \leq \sigma\norm{\nabla x}_H^2 + C_\sigma\norm x_*^2 \qquad\forall\,x\in V_1\,,
\eeq
and
\[
  \sp{\partial_t x(t)}{\mathcal Nx(t)}_{V_1}=\frac12\frac{d}{dt}\norm{\nabla\mathcal Nx(t)}_H^2
   \quad\text{for a.e.~}t\in(0,T)
\]
for every $x\in H^1(0,T; V^*)$ with $x_D=0$ almost everywhere in $(0,T)$.
We shall denote for simplicity 
\[
  H_0:=\{x\in H:\;x_D=0\}\,.
\]

The following assumptions on the data of the problem will be in force
throughout:
\begin{itemize}
  \item[(A1)] $\Psi\in C^2(\erre, \erre_+)$;
  \item[(A2)] there exist $c_1, c_2>0$ such that, for every $r\in\erre$,
  \[
  \Psi''(r)\geq -c_1\,,\qquad
  |\Psi''(r)|\leq c_2(1+|r|^2)\,,\qquad
  |\Psi'(r)|\leq c_2(1+\Psi(r))\,.
  \]
  \item[(A3)] $y_0\in L^{12}(\Omega,\cF_0; H)\cap L^6(\Omega,\cF_0; V)$ 
  and $\Psi(y_0)\in L^3(\Omega; L^1(D))$;
  \item[(A4)] $B:[0,T]\times H\to \cL^2(U,V)$ is  measurable
  and there exists a constant $L_B>0$ such that,
  for every $t\in[0,T]$,
  \begin{align*}
  \norm{B(t,x_1)-B(t,x_2)}_{\cL^2(U,H)}\leq L_B\norm{x_1-x_2}_H
  \quad&\forall\,x_1,x_2\in H\,,\\
  \norm{B(t,x)}_{\cL^2(U,H)}\leq L_B(1+\norm{x}_H)
  \quad&\forall\,x\in H\,,\\
  \norm{B(t,x)}_{\cL^2(U,V)}\leq L_B(1+\norm{x}_V)
  \quad&\forall\,x\in V\,.
  \end{align*}
  If $B$ is genuinely of multiplicative type, i.e.~if it is not constant 
  in the last variable, we further assume that the image of $B$ is contained
  in $\cL^2(U,H_0)$.
  \item[(A5)] for every $t\in[0,T]$, the operator
  $B(t,\cdot):H\to\cL^2(U,H)$ is of class~$C^1$.
\end{itemize}

\begin{remark}
  Let us comment on assumptions (A1)--(A5).
  In order for the state system \eqref{eq:1}--\eqref{eq:init} to be well-posed,
  one can require less stringent assumptions on the data (see for example
  \cite{scar-SCH, scar-SVCH}). However, in order to study the linearized
  and the adjoint problems, one needs some further regularity on
  the solution $y$ to the state equation, and for this reason (A1)--(A5) are in order.
  Let us point out that by the hypotheses (A1)--(A2) we can decompose $\Psi'$ as
  the sum of a continuous increasing function and a Lipschitz-continuous function
  as $\Psi'(r)=(\Psi'(r)+c_1r)-c_1r$, $r\in\erre$. Furthermore, note that (A3)--(A4) 
  are trivially satisfied for example when 
  $y_0\in V$ is nonnradom with $\Psi(y_0)\in L^1(D)$ and
  $B\in \cL^2(U,V)$ is time-independent
  and of additive type. 
  The reason why we assume 
  existence of higher moments on $y_0$ might not be 
  intuitive at this level and will be clarified later: let us mention that these assumptions
  will be needed to solve the linearized system and the adjoint problem,
  and that the hypothesis on the moment of order $6$ is ``optimal'' in this sense.
  In case of multiplicative noise,
  assumption (A4) corresponds to usual 
  boundedness and Lipschitz-continuity conditions on $B$,
  and the differentiability assumption (A5) is needed in order
  to analyse the linearized system. In particular,
  (A4)--(A5) imply that 
  \[
  \norm{DB(t,x)}_{\cL(H;\cL^2(U,H))}\leq L_B \quad\forall\,(t,x)\in [0,T]\times H\,.
  \]
\end{remark}

We define the set of admissible controls $\mathcal U$ as
\begin{align*}
  \mathcal U:=&\left\{u\in L^{12}_{\mathcal P}(\Omega; L^2(0,T; H))\cap
  L^6_{\mathcal P}(\Omega; L^2(0,T; V)):\right.\\
  &\left.\qquad\qquad\quad
  \norm{u}_{L^{12}(\Omega; L^2(0,T; H))\cap
  L^6(\Omega; L^2(0,T; V))}\leq C_0\right\}\,,
\end{align*}
where $C_0>0$ and $s\in(0,1/2)$ are fixed constants. 
It will be useful to introduce also the bigger set
\begin{align*}
  \mathcal U':=&\left\{u\in L^{12}_{\mathcal P}(\Omega; L^2(0,T; H))\cap
  L^6_{\mathcal P}(\Omega; L^2(0,T; V)):\right.\\
  &\left.\qquad\qquad\quad
  \norm{u}_{L^{12}(\Omega; L^2(0,T; H))\cap
  L^6(\Omega; L^2(0,T; V))}< 2C_0\right\}\,,
\end{align*}
which is open and bounded in
$L^{12}_{\mathcal P}(\Omega; L^2(0,T; H))\cap
L^6_{\mathcal P}(\Omega; L^2(0,T; V))$, and
$\mathcal U\subset \mathcal U'$.
Moreover, we define the cost functional 
\begin{align*}
  &J:L^2_{\mathcal P}(\Omega; C^0([0,T]; H))\times
  L^2_{\mathcal P}(\Omega; L^2(0,T; H))\to\erre_+\,,\\
  &J(y,u):=\frac{\alpha_1}{2}\E\int_Q|y-x_Q|^2 + 
  \frac{\alpha_2}{2}\E\int_D|y(T)-x_T|^2 +
  \frac{\alpha_3}{2}\E\int_Q|u|^2\,,
\end{align*}
where $\alpha_1, \alpha_2, \alpha_3\geq0$ are fixed constants and
\[
  \alpha_1x_Q \in L^6_{\mathcal P}(\Omega; L^6(0,T; H))\,, \qquad
  \alpha_2 x_T\in L^6(\Omega,\cF_T; V)\,.
\]

\begin{remark}
  The choices $\alpha_2 x_T\in L^6(\Omega; V)$ and 
  $\alpha_1 x_Q \in L^6(\Omega; L^6(0,T; H))$
  might look unnnatural to the reader at this level, 
  due to the form of the cost functional. However, this will be necessary in
  order to solve the adjoint system. Let us mention that conditions of this type are not new in
  literature of optimal control problems: see for example \cite{colli-sprek-optACDBC}
  for an analogous assumption in the context of the Allen-Cahn equation.
\end{remark}

As we have anticipated in Section~\ref{sec:intro}, we are interested in
minimizing $J(y,u)$ subject to the constraint $u\in\mathcal U$ 
and the state system \eqref{eq:1}--\eqref{eq:init}. We shall call optimal pair
any couple $(y,u)$ with $u\in \mathcal U$ satisfying \eqref{eq:1}--\eqref{eq:init}
and minimizing the cost functional $J$.

Under the hypotheses (A1)--(A4), we can prove that the state system is well-posed for every 
admissible control, and that the map $u\mapsto y$ is well-defined and Lipschitz-continuous.
These results are summarized in the following theorem.
\begin{theorem}
  \label{thm:1}
  Assume (A1)--(A4). Then for every $u\in \mathcal U'$ 
  there exists a unique pair $(y,w)$ with 
  \begin{gather}
    \label{state1}
    y \in L^{12}_{\mathcal P}\left(\Omega; C^0([0,T]; H)\cap  L^2(0,T; Z)\right)\,,\\
    \label{state1'}
    y \in L^6_{\mathcal P}\left(\Omega; L^\infty(0,T; V)\right)
    \cap L^3(\Omega; L^2(0,T; H^3(D)))\,,\\
    \label{state2}
    y-\int_0^\cdot B(s,y(s))\,dW(s) \in L^6(\Omega; H^1(0,T; V^*))\,,\\
    \label{state3}
    w \in L^3_{\mathcal P}(\Omega; L^2(0,T; V)))\,,\qquad
    \Psi'(y) \in 
    L^{3}_{\mathcal P}(\Omega; L^2(0,T; V))\,,
  \end{gather}
  such that $y(0)=y_0$ and, for every $\varphi\in V$, 
  for almost every $t\in(0,T)$, $\P$-almost surely,
  \begin{gather}
    \label{state5}
    \ip{\partial_t\left(y-\int_0^\cdot B(s,y(s))\,dW(s)\right)(t)}{\varphi}_V 
    + \int_D\nabla w(t)\cdot\nabla \varphi = 0\,,\\
    \label{state6}
    \int_Dw(t)\varphi = \int_D\nabla y(t)\cdot\nabla\varphi + \int_D\Psi'(y(t))\varphi - \int_Du(t)\varphi\,.
  \end{gather}
  Moreover, there exists a constant $M'>0$, only depending on
  $y_0$, $C_0$, $c_1$, $c_2$, $L_B$ and $Q$, such that,
  for every $u\in\mathcal U'$ and for any respective state $(y,w)$
  satisfying \eqref{state1}--\eqref{state6},
  \begin{align}
    \label{est1_state}
    \norm{y}_{L^{12}(\Omega; C^0([0,T]; H)\cap L^2(0,T; Z))}&\leq M'\,,\\
    \label{est1'_state}
    \norm{y}_{L^6(\Omega; L^\infty(0,T; V))\cap L^3(\Omega; L^2(0,T; H^3(D))}&\leq M'\,,\\
    \label{est2_state}
    \norm{w}_{L^3(\Omega; L^2(0,T; V))}+\norm{\Psi'(y)}_{L^3(\Omega; L^2(0,T; V))}
    &\leq M'\,.
  \end{align}
  Finally, there exists a constant $M>0$,
  only depending on
  $y_0$, $C_0$, $c_1$, $c_2$, $L_B$ and $Q$, such that, 
  for any $u_1,u_2\in \mathcal U'$
  and for any respective pairs $(y_1,w_1)$, $(y_2,w_2)$ 
  satisfying \eqref{state1}--\eqref{state6}, it holds
  \begin{align}
  \label{dep_cont}
  \norm{y_1-y_2}_{L^6(\Omega; C^0([0,T]; V^*))\cap L^6(\Omega; L^2(0,T; V))}
  &\leq M\norm{u_1-u_2}_{L^6(\Omega; L^2(0,T; V^*))}\,,\\
  \label{dep_cont'}
  \norm{y_1-y_2}_{L^2(\Omega; C^0([0,T]; H) \cap L^2(0,T; Z))}
  &\leq M \norm{u_1-u_2}_{L^6(\Omega; L^2(0,T; H))}\,.
  \end{align}
\end{theorem}

By Theorem~\ref{thm:1}, it is clear that uniqueness of $y$ holds for the
state system. Consequently, it is well-defined the control-to-state map
\[
  S: \mathcal U' \to
  L^6\left(\Omega; C^0([0,T]; H)\cap L^\infty(0,T; V)\cap L^2(0,T; Z)\right)\,,
\]
which is Lipschitz-continuous in the sense specified in \eqref{dep_cont}--\eqref{dep_cont'}. This allows us to 
introduce the reduced cost functional as
\[
  \tilde J : \mathcal U' \to \erre_+\,, \qquad
  \tilde J(u):=J(S(u),u)\,, \quad u\in \mathcal U'\,.
\]
The optimal control problem is thus equivalent to minimizing $\tilde J$ over 
$\mathcal U\subset \mathcal U'$.
The following definitions of optimal control are very natural.
\begin{definition}
  An optimal control is an element $u\in \mathcal U$ such that 
  \[
  \tilde J(u) \leq \tilde J(v) \qquad\forall\,v\in\mathcal U\,.
  \]
  A relaxed optimal control is a family
  \[
  ((\Omega^*, \cF^*, (\cF^*_t)_{t\in[0,T]},\P^*), W^*, x_Q^*, x_T^*, y_0^*, u^*, y^*, w^*)\,,
  \]
  where
  $(\Omega^*, \cF^*, (\cF^*_t)_{t\in[0,T]},\P^*)$ is
  a filtered probability space 
  satisfying the usual conditions, $W^*$ is
  a $(\cF^*_t)_t$-cylindrical Wiener process
  with values in $U$, $X_Q^*$ is a 
  $(\cF^*_t)_t$-progressively measurable 
  $L^2(0,T; H)$-valued process with the same law of $x_Q$,
  $x_T^*$ is a $\cF^*_T$-measurable $H$-valued random
  variable with the same law of $x_T$, $y_0^*$ is
  a $\cF^*_0$-measurable random variable 
  with the same law of $y_0$, $u^*$ is a process in
  the set $\mathcal U^*$ (defined as $\mathcal U$
  replacing $\Omega$ with $\Omega^*$),
  $(y^*,w^*)$ is the unique solution to the system
  \eqref{state1}--\eqref{state6} on $\Omega^*$
  with respect to the data $(W^*, y_0^*, u^*)$,
  and such that
  \[
  \tilde J^*(u^*):=\frac{\alpha_1}{2}\E{}^*\int_Q|y^*-x_Q^*|^2 + 
  \frac{\alpha_2}{2}\E{}^*\int_D|y^*(T)-x_T^*|^2
  +\frac{\alpha_3}{2}\E{}^*\int_Q|u^*|^2
   \leq \inf_{v\in\mathcal U}\tilde J(v)\,.
  \] 
\end{definition}
The first main result that we prove concerns with the existence of a relaxed optimal control.
Note however that existence of (strong) optimal controls is nontrivial, since the minimization problem
in not convex, hence uniqueness of optimal controls may fail. 
In case of uniqueness of optimal controls, existence of a strong optimal 
control can be proved for example by a well-known criterion on convergence in probability
due to Gy\"ongy--Krylov
(see \cite[Lem.~1.1]{gyo-kry}, \cite[Prop.~4.16]{hofma-parab}, and
\cite[Def.~2.4 and Thm.~2.5]{barbu-rock-contr}).
\begin{theorem}
  \label{thm:2}
  Assume (A1)--(A4). Then there exists a relaxed optimal control.
\end{theorem}

We focus now on the necessary conditions for optimality.
As we have anticipated, the first step consists in showing that $S$ is G\^ateaux-differentiable
in certain weak-sense and to characterize its weak derivative as the unique solution
of a linearized system. 
The following proposition ensures that the linearized system
is well-posed in a suitable variational sense.
\begin{proposition}
  \label{prop:lin}
  Assume (A1)--(A5). Then, for every $u\in \mathcal U'$ and for every
  $h\in L^6_{\mathcal P}(\Omega; L^2(0,T; H))$, setting $y:=S(u)$,
  there exists a unique pair $(z_h,\mu_h)$ with
  \begin{gather}
  \label{lin1}
  z_h \in L^2_{\mathcal P}\left(\Omega; 
  C^0([0,T]; H)\cap L^2(0,T; Z)\right)\,,\\
  \label{lin1'}
  z_h-\int_0^\cdot DB(s,y(s))z_h(s)\,dW(s) \in L^2(\Omega; H^1(0,T; Z^*))\,,\\
  \label{lin2}
  \mu_h \in L^2_{\mathcal P}(\Omega; L^2(0,T; H))\,,
  \end{gather}
  such that $z_h(0)=0$ and, for every $\varphi\in Z$, for almost every $t\in(0,T)$, $\P$-almost surely,
  \begin{gather}
    \label{lin3}
    \ip{\partial_t\left(z_h(t)-\int_0^t DB(s,y(s))z_h(s)\,dW(s)\right)}{\varphi}_Z 
    - \int_D\mu_h(t)\Delta\varphi = 0\,,\\
    \label{lin4}
    \int_D\mu_h(t)\varphi = \int_D\nabla z_h(t)\cdot\nabla\varphi
    +\int_D\Psi''(y(t))z_h(t)\varphi - \int_Dh(t)\varphi\,.
  \end{gather}
\end{proposition}

\begin{remark}
  Let us point out that \eqref{lin1}--\eqref{lin4} is the weak formulation of the linearized system,
  which can be obtained formally differentiating the state system \eqref{eq:1}--\eqref{eq:init}
  with respect to $u$ in the direction $h$, i.e.  
  \begin{align*}
    dz_h - \Delta \mu_h\,dt = DB(y)z_h\,dW \qquad&\text{in } (0,T)\times D\,,\\
    \mu_h=-\Delta z_h + \Psi''(y)z_h - h \qquad&\text{in } (0,T)\times D\,,\\
    \partial_{\bf n} z_h = \partial_{\bf n} \mu_h = 0 \qquad&\text{in } (0,T)\times\partial D\,,\\
     z_h(0)=0 \qquad&\text{in } D\,.
  \end{align*}
\end{remark}

We are now able to give a characterization of the weak G\^ateaux derivative of $S$
in terms of the unique solution to the linearized system \eqref{lin1}--\eqref{lin4}.
\begin{theorem}
  \label{thm:3}
  Assume (A1)--(A5). Then the control-to-state map $S$
  is weakly G\^ateaux-differentiable
  from $\mathcal U'$
  to $L^2(\Omega; C^0([0,T]; H))\cap L^2(\Omega; L^2(0,T; Z))$
  in the following sense: for every $u,h\in \mathcal U'$, as $\eps\searrow0$,
  \begin{align*}
  \frac{S(u+\eps h) - S(u)}{\eps}\to z_h \qquad&\text{in } 
  L^p(\Omega; L^2(0,T; V))\quad\forall\,p\in[1,2)\,,\\
  \frac{S(u+\eps h) - S(u)}{\eps}\wto z_h \qquad&\text{in } 
  L^2\left(\Omega; L^2(0,T; Z)\right)\,,\\
  \frac{S(u+\eps h)(t) - S(u)(t)}{\eps}\wto z_h(t) \qquad&\text{in } 
  L^2(\Omega; H) \quad\forall\,t\in[0,T]\,,
  \end{align*}
  where $z_h$ is the unique solution to the linearized system \eqref{lin1}--\eqref{lin4}.
\end{theorem}

The first natural necessary optimality condition is collected in the following result.
\begin{theorem}
  \label{thm:4}
  Assume (A1)--(A5), let $\bar u\in\mathcal U$ be an optimal control and $\bar y:=S(\bar u)$
  be the respective optimal state.
  Then 
  \[
  \alpha_1\E\int_Q(\bar y - x_Q)z_{v-\bar u} + \alpha_2\E\int_D(\bar y(T)-x_T)z_{v-\bar u}(T)
  +\alpha_3\E\int_Q\bar u(v-\bar u) \geq 0 \quad\forall\,v\in\mathcal U\,,
  \]
  where $z_{v-\bar u}$ is the unique solution to \eqref{lin1}--\eqref{lin4}
  with respect to the choice $h:=v-\bar u$.
\end{theorem}

The last result that we present is an alternative formulation of the 
first-order necessary conditions for optimality which does not involve
the solution $z$ to the linearized problem, but the solution
to the corresponding adjoint problem.
In this sense, the advantage is that the resulting
variational inequality that we obtain is much simpler to interpret.
The following proposition states that the adjoint problem
is well-posed in a suitable variational sense.
\begin{proposition}
  \label{prop:ad}
  Assume (A1)--(A5).
  Then for every $u\in \mathcal U'$, setting $y:=S(u)$,
  there exists a triple of processes $(p,\tilde p, q)$, with
  \begin{gather}
    \label{ad1}
    p\in C^0_w([0,T]; L^{2}(\Omega; V)) \cap 
    L^{2}_{\mathcal P}(\Omega; L^2(0,T; Z\cap H^3(D)))\,,\\
    \label{ad1'}
    \tilde p \in C^0_w([0,T]; L^{6}(\Omega; V^*)) \cap L^{6}_{\mathcal P}(\Omega; L^2(0,T; V))\,,\\
    \label{ad2}
    q\in L^2_{\mathcal P}(\Omega; L^2(0,T; \cL^2(U,V)))\,,
  \end{gather}
  such that, for every $\varphi\in Z$, $\P$-almost surely
  and for every $t\in[0,T]$,
  \begin{gather}
  \label{ad3}
  \int_D\tilde p(t)\varphi = \int_D\nabla p(t)\cdot\nabla\varphi\,,\\
  \label{ad4}
  \begin{split}
  &\int_Dp(t)\varphi -\int_t^T\!\!\int_D\tilde p(s)\Delta\varphi\,ds
  +\int_t^T\!\!\int_D\Psi''(y)\tilde p(s)\varphi\,ds\\
  &\qquad= \alpha_2\int_D(\bar y(T) - x_T)\varphi
  +\alpha_1\int_t^T\!\!\int_D(y-x_Q)(s)\varphi\,ds\\
  &\qquad+\int_t^T\left(DB(s,y(s))^*q(s), \varphi\right)_H\,ds
  -\int_D\left(\int_t^Tq(s)\,dW(s)\right)\varphi\,.
  \end{split}
  \end{gather}
  Moreover, 
  if $(p_1,\tilde p_1,q_1)$ and $(p_2,\tilde p_2, q_2)$
  are two solutions to \eqref{ad1}--\eqref{ad4}, then
  \[
  p_1-(p_1)_D = p_2-(p_2)_D\,, \qquad \tilde p_1 = \tilde p_2\,.
  \]
\end{proposition}

Our last result is a simplified version of the first-order necessary optimality conditions
which do not involve the solution $z$ to the linearized system, but the 
unique solution $\tilde p$ to the adjoint problem instead.

\begin{theorem}
\label{thm:5}
  Assume (A1)--(A5), let $\bar u\in\mathcal U$ be an optimal control and let $\bar y:=S(\bar u)$
  be the respective optimal state. Then the following variational inequality holds:
  \[
  \E\int_Q\left(\tilde p + \alpha_3 \bar u\right)(v-\bar u) \geq 0 \qquad\forall\,v\in\mathcal U\,,
  \]
  where $\tilde p$ is the unique second solution component satisfying \eqref{ad1}--\eqref{ad4}
  with respect to $(\bar u, \bar y)$.
  In particular, if $\alpha_3>0$, then $\bar u$ is the orthogonal projection of the point $-\frac{\tilde p}{\alpha_3}$
  on the closed convex set $\mathcal U$ in the Hilbert space $L^2(\Omega;L^2(0,T; H))$.
\end{theorem}

\begin{remark}
  The form of the cost functional $J$ considered in this paper is of standard quadratic tracking-type
  and is widely used in optimal control theory. However, let us point out that
  this choice is a particular case of the more general class of nonlinear performances
  \[
  J(y,u)=\E\int_0^TF_Q(t,y(t),u(t))\,dt + \E F_T(y(T))\,,
  \]
  where $F_Q:[0,T]\times H \times H\to \erre$ and 
  and $F_T: H\to \erre$ are
  $\cB([0,T])\otimes\cB(H)\otimes\cB(H)$--measurable
  and $\cB(H)$--measurable, respectively.
  The techniques used here can also be adapted to deal
  with such more general situations.
  For example,
  one can show existence 
  of relaxed optimal controls under very natural 
  lower semicontinuity assumptions on $F_Q$ and $F_T$.
  Furthermore, requiring that $F(t,\cdot,\cdot):H\times H\to\erre$
  and $F_T:H\to\erre$ are
  Fr\'echet-differentiable for every $t\in[0,T]$, with 
  \begin{align*}
  \norm{D_yF_Q(t,y,u)}_H + \norm{D_uF_Q(t,y,u)}_H&\lesssim 1 + \norm{y}_H + \norm{u}_H\\
  \norm{D_yF_T(y)}_H&\lesssim 1 + \norm{y}_H
  \end{align*}
  for every $(t,y,u)\in[0,T]\times H\times H$, necessary conditions 
  for optimality can be also studied.
  Note however that in the case of nonlinear performance,
  the resulting variational inequality in Theorem~\ref{thm:5}
  would not give a characterization of the optimal controls in terms 
  of orthogonal projection on $\mathcal U$.
\end{remark}

\section{Well-posedness of the state system}
\label{sec:wp}
\setcounter{equation}{0}

This section is devoted to the proof of Theorem~\ref{thm:1},
ensuring the the state system is well-posed.

For any $\lambda>0$, we consider the approximated problem
\[
  \begin{cases}
  dy_\lambda - \Delta w_\lambda\,dt = B(y_\lambda)\,dW \quad&\text{in } (0,T)\times D\,,\\
  w_\lambda=-\Delta y_\lambda + \Psi'_\lambda(y_\lambda) - u \quad&\text{in } (0,T)\times D\,,\\
  \partial_{\bf n}y_\lambda=
  \partial_{\bf n}w_\lambda=0 \quad&\text{in } (0,T)\times \partial D\,,\\
  y_\lambda(0)=y_0 \quad&\text{in } D\,,
  \end{cases}
\]
where $\Psi'_\lambda$ is a Lipschitz-continuous 
smooth Yosida-type approximation of $\Psi'$. 
The classical variational theory ensures the existence and uniqueness
of an approximated solution $y_\lambda \in L^{12}_{\mathcal P}(\Omega; C^0([0,T];H)\cap L^2(0,T;Z))$.
Arguing as in \cite{scar-SCH, scar-SVCH}, It\^o's formula for
the square of the $H$-norm and the linear growth 
condition on $B$ in assumption (A4) yields, together with the Gronwall lemma, that
\beq\label{est_wp1}
  \norm{y_\lambda}_{L^{12}(\Omega; C^0([0,T]; H)\cap L^2(0,T; Z))}\leq M'\,,
\eeq
where the constant $M'>0$ only depends on
$y_0$, $C_0$, $c_1$, $c_2$, $L_B$ and $Q$.
Furthermore, writing 
It\^o's formula for the free-energy functional (see again \cite{scar-SCH}) yields
\beq\label{ito_aux}
\begin{split}
&\frac12\int_D|\nabla y_\lambda(t)|^2 + \int_D\Psi_\lambda(y_\lambda(t))
+\int_{Q_t}\nabla w_\lambda\cdot\nabla(w_\lambda + u)\\
&= \frac12\int_D|\nabla y_0|^2 + \int_D\Psi_\lambda(y_0)
+\int_0^t\left((w_\lambda + u)(s), B(s, y_\lambda(s))\right)_H\,dW(s)\\
&\qquad+\frac12\int_0^t\norm{\nabla B(s,y_\lambda(s))}^2_{\cL^2(U,H)}\,ds
+\frac12\sum_{k=0}^\infty\int_{Q_t}\Psi_\lambda''(y_\lambda)|B(\cdot, y_\lambda)e_k|^2\,.
\end{split}
\eeq
We want now to take power $3$ at both sides, and then supremum in time
and expectations.
Note that the trace term on the right-hand side can be estimated
thanks to the H\"older inequality, the Sobolev embedding $V\embed L^6(D)$
and the growth conditions (A2) and (A4) on $\Psi''$ and $B$, as
\begin{align*}
  &\sum_{k=0}^\infty\int_{Q_t}\Psi_\lambda''(y_\lambda)|B(\cdot, y_\lambda)e_k|^2\\
  &\lesssim
  \norm{B(\cdot, y_\lambda)}^2_{L^2(0,t; \cL^2(U,H))} + 
  \int_0^t\norm{y_\lambda(s)}_{L^4(D)}^2\norm{B(\cdot, y_\lambda)}^2_{\cL^2(U,V)}\,ds\\
  &\lesssim 1 + \norm{y_\lambda}^2_{L^2(0,T; V)} 
  + \norm{y_\lambda}^4_{L^4(0,T; V)}\,.
\end{align*}
Since by interpolation we have 
\[
\norm{y_\lambda}_{L^4(0,t; V)}\lesssim\norm{y_\lambda}_{L^\infty(0,T; H)\cap L^2(0,T; Z)}\,,
\]
the right-hand side is uniformly bounded in $L^3(\Omega)$
by the estimate \eqref{est_wp1}.
Furthermore, for the stochastic integral we note that 
\[
  (w_\lambda + u, B)_H = (w_\lambda-(w_\lambda)_D, B)_H + |D|(w_\lambda)_DB_D + (u,B)\,.
\]
The Burkholder-Davis-Gundy, Young, and
Poincar\'e-Wirtinger inequalities imply, together with assumption (A4), that 
\begin{align*}
  &\E\sup_{r\in[0,t]}\left|\int_0^r\left((w_\lambda + u)(s),
  B(s, y_\lambda(s))\right)_H\,dW(s)\right|^3\\
  &\lesssim\E\left(\int_0^t\left(\norm{(w_\lambda-(w_\lambda)_D)(s)}_H^2+\norm{u(s)}_H^2\right)
  \norm{B(s,y_\lambda(s))}_{\cL^2(U,H)}^2\,ds\right)^{3/2}\\
  &\qquad+\E\left(\int_0^t|(w_\lambda(s))_D|^2\norm{B_D(s,y_\lambda(s))}_{\cL^2(U,\erre)}^2\,ds\right)^{3/2}\,.
\end{align*}
Now, in case of multiplicative noise we have that $B_D=0$ by (A4) and
the second term on the right-hand side vanishes, 
so that we can continue the estimate by
\begin{align*}
  &\E\left[\norm{B(\cdot,y_\lambda)}_{L^\infty(0,T; \cL^2(U,H))}^3
  \left(\norm{w_\lambda-(w_\lambda)_D}_{L^2(0,t;H)}^3
  +\norm{u}_{L^2(0,T; H)}^3\right)\right]\\
  &\lesssim 1+\delta\norm{\nabla w_\lambda}^6_{L^6(\Omega; L^2(0,T; H))}
  +C_\delta\norm{y_\lambda}^6_{L^6(\Omega; C^0([0,T]; H))}
  +\norm{u}^6_{L^6(\Omega; L^2(0,T; H))}\,.
\end{align*}
In case of additive noise, on the right-hand side we obtain the further contribution 
\[
  \norm{B}^3_{L^\infty(0,T; \cL^2(U,V^*))}\norm{(w_\lambda)_D}_{L^3(\Omega; L^2(0,t))}^3\lesssim
  t^{3/2}\norm{(w_\lambda)_D}^3_{L^3(\Omega; L^\infty(0,t))}\,.
\]
We go back now to \eqref{ito_aux}, and take power $3$, supremum in $t\in[0,T_0]$
for a certain $T_0\in(0,T)$
and expectations.
Since $w_\lambda=-\Delta y_\lambda + \Psi'_\lambda(y_\lambda) - u$, by (A2)
the term involving $\Psi_\lambda(y_\lambda)$ on the left-hand side of \eqref{ito_aux}
yields a bound from below for $(w_\lambda)_D$ in $L^3(\Omega; L^\infty(0,T_0))$: hence,
choosing $\delta>0$ and $T_0$ small enough,
rearranging the terms, and using the Gronwall lemma yields
\begin{align*}
&\norm{y_\lambda}_{L^6(\Omega; L^\infty(0,T_0; V))}^6
+\norm{\Psi_\lambda(y_\lambda)}_{L^3(\Omega; L^\infty(0,T_0; L^1(D)))}
+\norm{\nabla w_\lambda}_{L^6(\Omega; L^2(0,T_0; H))}^6\\
&\lesssim 1 + \norm{y_\lambda}^{12}_{L^{12}(\Omega; C^0([0,T]; H)\cap L^2(0,T; Z))}
+\norm{u}_{L^6(\Omega; L^2(0,T; V))}^{6}\,,
\end{align*}
where the implicit constant is independent of $\lambda$.
Taking these remarks into account, noting that 
$\norm{w_\lambda}_V\lesssim\norm{\nabla w_\lambda}_H + |(w_\lambda)_D|$,
using \eqref{est_wp1} and the fact that $u\in\mathcal U'$, we infer by
using a classical patching-in-time technique that 
\beq
  \label{est_wp2}
  \norm{y_\lambda}_{L^6(\Omega; L^\infty(0,T; V))} 
  + \norm{w_\lambda}_{L^3(\Omega; L^2(0,T; V))}\leq M'\,.
\eeq
Now, by comparison in the equation for the chemical potential we 
deduce that 
\[
  \norm{\Psi'_\lambda(y_\lambda)}_{L^3(\Omega; L^2(0,T; H))}\leq M'\,,
\]
while by (A2), the embedding $V\embed L^6(D)$, and interpolation
we have 
\[
  \norm{\nabla \Psi'(y_\lambda)}_{L^2(0,T; H)}=
  \norm{\Psi''_\lambda(y_\lambda)\nabla y_\lambda}\lesssim
  1 + \norm{y_\lambda}_{L^\infty(0,T; V)}\norm{y_\lambda}^2_{L^\infty(0,T; H)\cap L^2(0,T; Z)}\,,
\]
so that by \eqref{est_wp1}--\eqref{est_wp2}
\[
\norm{\Psi'_\lambda(y_\lambda)}_{L^3(\Omega; L^2(0,T; V))}\leq M'\,.
\]
By elliptic regularity we infer then also 
\beq
  \label{est_wp3}
  \norm{\Psi'(y_\lambda)}_{L^3(\Omega; L^2(0,T; V))}+
  \norm{y_\lambda}_{L^3(\Omega; L^2(0,T; H^3(D)))}\leq M'\,.
\eeq
It is now a standard matter to pass to the limit as $\lambda\searrow0$
in the approximated problem, and recover \eqref{state1}--\eqref{est2_state}:
for further details we refer to \cite{scar-SCH, scar-SVCH}.

It only remains to prove the continuous dependence property \eqref{dep_cont'},
as \eqref{dep_cont} has already been proved in \cite{scar-SVCH}.
To this end, note that 
  \begin{align*}
  &d(y_1-y_2)-\Delta(w_1-w_2)\,dt =(B(t,y_1)-B(t,y_2))\,dW\,,\\
  &w_1-w_1=-\Delta(y_1-y_2) 
  +\Psi'(y_1)-\Psi'(y_2)-(u_1-u_2)\,,
  \end{align*}
  so that It\^o's formula for the square of the $H$-norm yields
  \begin{align*}
  &\frac12\norm{(y_1-y_2)(t)}_H^2 + \int_{Q_t}|\Delta(y_1-y_2)|^2\\
  &=\int_{Q_t}\left(\Psi'(y_1)-\Psi'(y_2)\right)\Delta(y_1-y_2)
  -\int_{Q_t}(u_1-u_2)\Delta(y_1-y_2)\\
  &\qquad+\frac12\norm{B(\cdot,y_1)-B(\cdot,y_2)}_{L^2(0,t;\cL^2(U,H))}^2\\
  &\qquad+\int_0^t\left((y_1-y_2)(s), B(s,y_1(s))-B(s,y_2(s))\right)_H\,dW(s)\,.
  \end{align*}
  Using the Burkholder-Davis-Gundy and Young inequalities, and employing
  the Lipschitz continuity of $B$, we have
  \begin{align*}
  &\norm{y_1-y_2}_{L^2(\Omega; C^0([0,t]; H))}^2 
  + \norm{\Delta(y_1-y_2)}_{L^2(\Omega; L^2(0,t; H))}^2\\
  &\lesssim\E\int_Q|u_1-u_2|^2+
  \E\int_Q|\Psi'(y_1)-\Psi'(y_2)|^2
  +\norm{y_1-y_2}^2_{L^2(\Omega; L^2(0,t; H))}\\
  &\qquad+\delta\norm{y_1-y_2}_{L^2(\Omega; C^0([0,t]; H))}^2
  +C_\delta \norm{y_1-y_2}^2_{L^2(\Omega; L^2(0,t; H))}\,,
  \end{align*}
  for every $\delta >0$ and a suitable constant $C_\delta >0$.
  By the mean-value theorem, assumption (A2), the H\"older inequality and
  the embedding $V\embed L^6(D)$,
  \begin{align*}
  &\E\int_Q|\Psi'(y_1)-\Psi'(y_2)|^2\lesssim\E\int_Q\left(1+|y_1|^4 + |y_2|^4\right)|y_1-y_2|^2\\
  &\qquad\lesssim\E\int_0^T\left(1+\norm{y_1(s)}^4_{L^6(D)} + \norm{y_2(s)}^4_{L^6(D)}\right)
  \norm{(y_1-y_2)(s)}^2_{L^6(D)}\,ds\\
  &\qquad\lesssim\E\left(1+\norm{y_1}_{L^\infty(0,T; V)}^4 + \norm{y_2}^4_{L^\infty(0,T; V)}\right)
  \norm{y_1-y_2}^2_{L^2(0,T; V)}\\
  &\qquad\lesssim\left(1+\norm{y_1}_{L^6(\Omega; L^\infty(0,T; V))}^4 
  + \norm{y_2}^4_{L^6(\Omega; L^\infty(0,T; V))}\right)
  \norm{y_1-y_2}^2_{L^6(\Omega; L^2(0,T; V))}\,.
  \end{align*}
  Hence, \eqref{dep_cont'}
  follows rearranging the terms and using \eqref{dep_cont} and \eqref{est1'_state}.

\section{Existence of a relaxed optimal control}
\label{sec:exist}
\setcounter{equation}{0}

This section is devoted to the proof of Theorem~\ref{thm:2}:
we show that a relaxed optimal control always exists. 

Let $(u_n)_n\subset\mathcal U$ be a minimizing sequence for the reduced cost functional $\tilde J$,
and set $(y_n, w_n)$ as the respective solution to \eqref{state1}--\eqref{state6}.
Then, by definition of $\mathcal U$ and the uniform estimates \eqref{est1_state}--\eqref{est2_state}, 
we deduce that there exists a positive constant $c$, independent of $n$, such that 
\begin{gather*}
    \norm{u_n}_{L^{12}(\Omega; L^2(0,T; H))\cap 
    L^6(\Omega; L^2(0,T; V))}\leq C_0\,,\\
    \norm{y_n}_{L^{12}(\Omega; C^0([0,T]; H)\cap L^2(0,T; Z))\cap
    L^6(\Omega; L^\infty(0,T; V))\cap
    L^3(\Omega; L^2(0,T; H^3(D)))}\leq c\,,\\
    \norm{w_n}_{L^3(\Omega; L^2(0,T; V))} +\norm{\Psi'(y_n)}_{L^3(\Omega; L^2(0,T; V))}\leq c\,.
\end{gather*}
Recalling hypothesis (A4), we also deduce that 
\[
  \norm{B(\cdot, y_n)}_{L^6(\Omega; L^\infty(0,T; \cL^2(U,V)))}\leq c\,.
\] 
Hence, by \cite[Lem.~2.1]{fland-gat}, for every $s\in(0,1/2)$, there
exists $c_s>0$, independent of $n$, such that
\[
  \norm{B(\cdot, y_n)\cdot W}_{
  L^6(\Omega; W^{s,6}(0,T; V))}\leq c_s\,,
\]
where we have used the classical notation $\cdot$ for the stochastic integral.
Since $1-\frac12>s-\frac16$, we have that 
$H^1(0,T; V^*)\embed W^{s,6}(0,T; V^*)$ by the Sobolev embeddings, and
by comparison in \eqref{state5} we infer that 
\[
  \norm{y_n}_{L^6(\Omega; W^{s,6}(0,T;V^*))}\leq c_s\,.
\]
Let us define now $\pi_{y_n}$ as the law of $y_n$ on $C^0([0,T]; H)\cap L^2(0,T; Z)$ 
and show that $(\pi_{y_n})_n$ is tight.
Fixing now $s\in(1/6, 1/2)$ so that $6s>1$,
by \cite[Sec.~8, Cor.~4--5]{simon}
we have the compact inclusions
\begin{align*}
  L^2(0,T; H^3(D)\cap Z)\cap W^{s,6}(0,T; V^*)&\cembed L^2(0,T; Z)\,, \\ 
  L^\infty(0,T; V)\cap W^{s,6}(0,T; V^*)&\cembed C^0([0,T]; H)\,.
\end{align*}
If we define the space
\[
  \mathcal W:= L^\infty(0,T; V)\cap L^2(0,T; H^3(D)\cap Z)\cap W^{s,6}(0,T; V^*)\,,
\]
we deduce that $\mathcal W\cembed C^0([0,T]; H)\cap L^2(0,T; Z)$ compactly and also the estimate
\[
  \norm{y_n}_{L^3(\Omega;\mathcal W)}\leq c\,.
\]
This ensures by a standard argument that 
$(\pi_{y_n})_n$ is tight on $C^0([0,T]; H)\cap L^2(0,T; Z)$. 
Indeed, if $B_R$ denotes the closed ball of radius $R>0$
in $\mathcal W$, for any $R>0$, we have that $B_R$ is compact in 
$C^0([0,T]; H)\cap L^2(0,T; Z)$, and by Markov's inequality
\[
  \pi_{y_n}(B_R^c)=\P\left\{\norm{y_n}_{\mathcal W}^3>R^3\right\}\leq\frac1{R^3}
  \norm{y_n}_{L^3(\Omega;\mathcal W)}^3\leq \frac{c^3}{R^3} \qquad\forall\,n\in\enne\,,
\]
from which the tightness of $(\pi_{y_n})_n$. Similarly, by \cite[Sec.~8, Cor.~4--5]{simon}
we also have the compact inclusion
$W^{s,6}(0,T; V) \cembed C^0([0,T]; H)$,
so that an entirely analogous argument yields that the laws
of $(B(\cdot, y_n)\cdot W)_n$ 
on the space $C^0([0,T]; H)$ are tight. Moreover,
denoting by $L^2_w(0,T; V)$ the space $L^2(0,T; V)$ endowed with 
its weak topology, it is clear that the laws of $(u_n)_n$
on $L^2_w(0,T; V)$ are tight.

Now, taking into account the remarks above, we deduce in particular that the septuple
$(y_n, u_n, B(\cdot, y_n)\cdot W, W, y_0, x_Q, x_T)_n$ is tight on the space
\[
  C^0([0,T]; H) \times
  L^2_w(0,T; V) \times C^0([0,T]; H) \times C^0([0,T]; U)
  \times V \times L^2(0,T; H) \times H\,.
\]
By Skorokhod theorem (see \cite[Thm.~2.7]{ike-wata} 
and \cite[Thm.~1.10.4, Add.~1.10.5]{vaa-well})
and the Jakubowski-Skorokhod version
(see \cite[Thm.~2.7.1]{fei-hof}),
there is a probability space $(\Omega^*, \cF^*, \P^*)$,
a sequence of maps $(\phi_n)_n$,
where $\phi_n:(\Omega^*,\cF^*)\to
(\Omega,\cF)$ are measurable and satisfy $\P= \P^*\circ\phi_n^{-1}$ for every $n\in\enne$,
and measurable random variables $(y^*,u^*, I^*, W^*, y_0^*, x_Q^*, x_T^*)$
defined on $(\Omega^*, \cF^*)$ with values in
\[
  C^0([0,T]; H) \times
  L^2(0,T; V) \times C^0([0,T]; H) \times C^0([0,T]; U)
  \times V \times L^2(0,T; H) \times H\,,
\]
such that 
\begin{align*}
  y_n^*:=y_n\circ\phi_n \to y^* \qquad&\text{in } C^0([0,T]; H) \quad\P^*\text{-a.s.}\,,\\
  u_n^*:=u_n\circ\phi_n \wto u^* \qquad&\text{in } L^2(0,T; V) \quad\P^*\text{-a.s.}\,,\\
  I_n^*:=(B(\cdot, y_n)\cdot W)\circ\phi_n \to I^* 
  \qquad&\text{in } C^0([0,T]; H) \quad\P^*\text{-a.s.}\,,\\
  W_n^*:=W\circ\phi_n \to W^* \qquad&\text{in } C^0([0,T]; U)\quad\P^*\text{-a.s.}\,,\\
  y_{0,n}^*:=y_0\circ\phi_n \to y_0^* \qquad&\text{in } V\quad\P^*\text{-a.s.}\,,\\
  x_{Q,n}^*:=x_Q\circ\phi_n \to x_Q^* \qquad&\text{in } L^2(0,T; H)\quad\P^*\text{-a.s.}\,,\\
  x_{T,n}^*:=x_T\circ\phi_n \to x_T^* \qquad&\text{in } H\quad\P^*\text{-a.s.}\,.
\end{align*}
Since the sequence $(y_0, x_Q, x_T)_n$ is constant, it is clear that 
the laws of $(y_0^*, x_Q^*, x_T^*)$ and $(y_0,x_Q,x_T)$ coincide.
Moreover, 
setting $w_n^*:=w_n\circ w_n$, 
since the maps $(\phi_n)_n$ preserve the laws, we readily deduce that 
\begin{gather*}
    \norm{y_n^*}_{L^{12}(\Omega^*; C^0([0,T]; H)\cap L^2(0,T; Z))\cap
    L^6(\Omega^*; L^\infty(0,T; V))\cap
    L^3(\Omega^*; L^2(0,T; H^3(D)))}\leq c\,,\\
    \norm{u_n^*}_{L^{12}(\Omega^*; L^2(0,T; H))\cap 
    L^6(\Omega^*; L^2(0,T; V))}\leq C_0\,,\\
    \norm{I_n^*}_{L^6(\Omega^*; W^{s,6}(0,T; V))}\leq c_s\,,\\
    \norm{w_n^*}_{L^3(\Omega^*; L^2(0,T; V))} +\norm{\Psi'(y_n^*)}_{L^3(\Omega^*; L^2(0,T; V))}\leq c\,,
\end{gather*}
hence in particular that 
\begin{gather*}
    y^* \in L^{12}(\Omega^*; C^0([0,T]; H)\cap  L^2(0,T; Z))\,,\\
    y^* \in L^6(\Omega^*;L^\infty(0,T; V)
    \cap L^3(\Omega^*; L^2(0,T; H^3(D))\,,\\
    u^* \in \mathcal U^*\,,\\
    I^* \in L^6(\Omega^*; W^{s,6}(0,T; V))
\end{gather*}
and
\begin{align*}
  y_n^*\to y^* \qquad&\text{in } 
  L^p(\Omega^*; C^0([0,T]; H))^2 \quad\forall\,p\in[1,12)\,,\\
  y_n^*\wto y^* \qquad&\text{in } L^{12}(\Omega^*; L^2(0,T; Z))
  \cap L^3(\Omega^*; L^2(0,T; H^3(D))) \,,\\
  u_n^* \wto u^* \qquad&\text{in } L^{12}(\Omega^*; L^2(0,T; H))\cap L^6(\Omega^*; L^2(0,T; V))\,,\\
  w_n^*\wto w^*\qquad&\text{in } L^6(\Omega^*; L^2(0,T; V))\,,\\
  \Psi'(y_n^*) \wto \xi^* \qquad&\text{in } L^{3}(\Omega^*; L^2(0,T; V))\,,
\end{align*}
for some 
\[
  w^* \in L^6(\Omega^*; L^2(0,T; V)\,,\qquad
    \xi^* \in L^{3}(\Omega^*; L^2(0,T; V))\,.
\]
The strong-weak closure of the maximal monotone operator
$r\mapsto\Psi'(r)+c_1r$, $r\in\erre$, ensures that $\xi^*=\Psi'(y^*)$
almost everywhere. Moreover, by (A4) we have that 
\[
  B(\cdot, y_n^*) \to B(\cdot, y^*) \qquad\text{in }
  L^p(\Omega^*; L^2(0,T; \cL^2(U,H))) \quad\forall\,p\in[1,12)\,.
\]
Now, defining the filtrations $(\cF^*_{n,t})_{t\in[0,T]}$ and 
$(\cF^*_t)_{t\in[0,T]}$ as
\[
  \cF^*_{n,t}:=\sigma(W^*_n(s))_{s\in[0,t]}\,, \quad
  \cF^*_t:=\sigma(y^*(s), u^*(s), I^*(s), W^*(s))_{s\in[0,T]}\,, \qquad t\in[0,T]\,,
\]
using classical representation theorems for martingales
(see for example the arguments in \cite[\S~4]{vall-zimm}) it is possible to show
that $W^*_n$ is a $(\cF^*_{n,t})_t$-cylindrical Wiener process on $U$, 
$W^*$ is a $(\cF^*_{t})_t$-cylindrical Wiener process on $U$, and that 
\[
  I_n^*=\int_0^\cdot B(s,y_n^*(s))\,dW^*_n(s)\,, \qquad I^*=\int_0^\cdot B(s,y^*(s))\,dW^*(s)\,.
\]
Since $(y_n^*, w_n^*)$ satisfies the variational 
formulation \eqref{state5}--\eqref{state6} on the space $\Omega^*$
with respect to $(y_{0,n}^*, u_n^*)$, 
passing to the weak limit
it follows that $(y^*,w^*)$ is the unique solution to \eqref{state1}--\eqref{state6}
on the probability space $(\Omega^*,\cF^*, \P^*)$
corresponding to $(y_0^*, u^*)$.
Using the weak lower semicontinuity of $J$,
the fact that $\phi_n$ preserves the law for every $n$,
and the definition of 
the minimizing sequence $(u_n)_n$, we deduce that
\begin{align*}
  &\tilde J^*(u^*)=
  \frac{\alpha_1}{2}\E{}^*\int_Q|y^*-x_Q^*|^2 + 
  \frac{\alpha_2}{2}\E{}^*\int_D|y^*(T)-x_T^*|^2 +
  \frac{\alpha_3}{2}\E{}^*\int_Q|u^*|^2\\
  &\leq\liminf_{n\to\infty}
  \frac{\alpha_1}{2}\E{}^*\int_Q|y_n^*-x_{Q,n}^*|^2 + 
  \frac{\alpha_2}{2}\E{}^*\int_D|y_n^*(T)-x_{T,n}^*|^2 +
  \frac{\alpha_3}{2}\E{}^*\int_Q|u_n^*|^2\\
  &=\liminf_{n\to\infty}\tilde J(u_n)=
  \lim_{n\to\infty}\tilde J(u_n)=\inf_{v\in\mathcal U}\tilde J(v)\,,
\end{align*}
so that $u^*$ is a relaxed optimal control.

\section{The control-to-state map}
\label{sec:map}
\setcounter{equation}{0}

In this section we study the G\^ateaux differentiability of the control-to-state map
and we prove the first version of first-order necessary conditions for optimality.

\subsection{Existence-uniqueness of the linearized system}
\label{subsec:lin}
We prove here Proposition~\ref{prop:lin}. 
Let $u\in \mathcal U'$ be given, set 
$y:=S(u)$, and let
$h\in L^6_{\mathcal P}(\Omega; L^2(0,T; H))$. We show that the linearized system
\eqref{lin1}--\eqref{lin2} admits a unique solution $z_h$.

\noindent{\bf Uniqueness.} 
For $i=1,2$, let 
\[
(z_h^i,\mu_h^i)\in L^2_{\mathcal P}\left(\Omega;C^0([0,T]; H)\cap L^2(0,T;Z)\right) 
\times L^2_{\mathcal P}(\Omega; L^2(0,T; H))\,,
\]
such that $(z_h^i,\mu_h^i)$ satisfy \eqref{lin3}--\eqref{lin4}. Then we have, 
in the variational sense in the triple $(Z,H,Z^*)$,
\begin{align*}
    d(z^1_h-z_h^2) - \Delta (\mu^1_h-\mu_h^2)\,dt 
    = DB(y)\left(z^1_h-z_h^2\right)\,dW \qquad&\text{in } (0,T)\times D\,,\\
    \mu_h^1-\mu_h^2=-\Delta (z^1_h-z_h^2) 
    + \Psi''(y)(z_h^1-z_h^2) \qquad&\text{in } (0,T)\times D\,,\\
    \partial_{\bf n} (z^1_h-z_h^2) = \partial_{\bf n} (\mu^1_h-\mu_h^2) = 0 
    \qquad&\text{in } (0,T)\times\partial D\,,\\
     (z^1_h-z_h^2)(0)=0 \qquad&\text{in } D\,.
\end{align*}
Integrating on $D$ the first equation, it follows from (A4) that 
$(z_h^1-z_h^2)_D=0$. Hence, It\^o's formula for the 
square of the $V^*$-norm yields
\begin{align*}
  &\frac12\norm{\nabla\mathcal N(z^1_h-z_h^2)(t)}_H^2 + \int_{Q_t}|\nabla(z^1_h-z_h^2)|^2
  +\int_{Q_t}\Psi''(y)|(z^1_h-z_h^2)|^2\\
  &=\frac12\int_0^t\operatorname{Tr}\left(DB(s,y(s))\left(z^1_h-z_h^2\right)(s)^*\circ
  \mathcal N\circ DB(s,y(s))\left(z^1_h-z_h^2\right)(s)\right)\,ds\\
  &+\int_0^t\left(\mathcal N(z_h^1-z_h^2)(s),DB(s,y(s))\left(z^1_h-z_h^2\right)(s)\right)_H\,dW(s)\,.
\end{align*}
Now, by the uniform boundedness of $DB$,
the first term on the right-hand side is bounded by 
\begin{align*}
  &\frac12\norm{DB(\cdot, y)(z_h^1-z_h^2)}_{L^2(0,t; \cL^2(U,V^*))}^2\lesssim_{L_B}
  1 + \norm{z_h^1-z_h^2}_{L^2(0,t; H)}^2\\
  &\qquad\leq
  \delta\norm{\nabla(z_h^1-z_h^2)}_{L^2(0,t; H)}^2
  + C_\delta\norm{z_h^1-z_h^2}_{L^2(0,t; V^*)}^2
\end{align*}
for every $\delta >0$ and a certain $C_\delta >0$,
while the second term on the right-hand side
can be estimated using the Burkholder-Davis-Gundy and Young
inequalities as
\begin{align*}
  &\E\sup_{r\in[0,t]}\left|
  \int_0^r\left(\mathcal N(z_h^1-z_h^2)(s),DB(s,y(s))
  \left(z^1_h-z_h^2\right)(s)\right)_H\,dW(s)\right|\\
  &\lesssim
  \delta\norm{z_h^1-z_h^2}_{L^2(\Omega; C^0([0,t]; V^*))}^2
  +C_\delta\norm{DB(\cdot,y)
  \left(z^1_h-z_h^2\right)}^2_{L^2(\Omega; L^2(0,t; \cL^2(U,H)))}\\
  &\lesssim
  \delta\norm{z_h^1-z_h^2}_{L^2(\Omega; C^0([0,t]; V^*))}^2
  +\delta\norm{\nabla(z_h^1-z_h^2)}_{L^2(\Omega; L^2(0,t; H))}^2\\
  &\qquad+C_\delta\norm{z_h^1-z_h^2}_{L^2(\Omega; L^2(0,t; V^*))}^2\,.
\end{align*}
Furthermore, since $\Psi''\geq -c_1$, choosing $\delta$ sufficiently small and
rearranging the terms
yields, the Young inequality,
\begin{align*}
  \frac12\norm{z^1_h-z_h^2}^2_{L^2(\Omega; C^0([0,t]; H))} &+ 
  \E\int_{Q_t}|\nabla(z^1_h-z_h^2)|^2
  \lesssim c_1\E\int_{Q_t}|z^1_h-z_h^2|^2\\
  &\lesssim\sigma\E\int_{Q_t}|\nabla(z^1_h-z_h^2)|^2
  +\tilde c_\sigma\E\int_0^t\norm{\nabla\mathcal N(z^1_h-z_h^2)(s)}_H^2\,ds
\end{align*}
for every $\sigma>0$ and a certain $\tilde c_\sigma>0$. 
Taking $\sigma$ small enough, 
the Gronwall lemma yields then $z_h^1=z_h^2$, hence also
$\mu_h^1=\mu_h^2$ by comparison in the system, from which uniqueness.

\noindent{\bf Approximation.}
Let us focus on existence.
To this end, we consider the approximation
\begin{align*}
    dz^n_h - \Delta \mu^n_h\,dt = DB(y)z_h^n\,dW \qquad&\text{in } (0,T)\times D\,,\\
    \mu_h^n=-\Delta z^n_h + \Psi_n''(y)z_h - h \qquad&\text{in } (0,T)\times D\,,\\
    \partial_{\bf n} z^n_h = \partial_{\bf n} \mu_h = 0 \qquad&\text{in } (0,T)\times\partial D\,,\\
     z^n_h(0)=0 \qquad&\text{in } D\,,
\end{align*}
where $\Psi''_n:=T_n\circ\Psi''$ and $T_n:\erre\to\erre$ is the truncation operator at level $n$. i.e.
\[
  T_n(r):=\begin{cases}
  n \quad&\text{if } r>n\,,\\
  r &\text{if } |r|\leq n\,,\\
  -n &\text{if } r<-n\,,
  \end{cases}
  \qquad r\in\erre\,.
\]
Since $\Psi''_n(y)\in L^\infty(\Omega\times Q)$, it is not difficult to check that such
approximated problem admits a unique solution $(z_h^n,\mu_h^n)$ satisfying \eqref{lin1}--\eqref{lin4}
with $\Psi''_n$ instead of $\Psi''$. Indeed, one can reformulate the problem in the 
variational triple $(Z,H,Z^*)$ as
\[
  dz_h^n + A_n z_h^n\,dt = DB(t,y)z_h^n\,dW\,, \qquad z_h^n(0)=0\,,
\]
where $A_n:\Omega\times[0,T]\times Z\to Z^*$ is given by 
\[
  \ip{A_n(\omega,t,x)}{\varphi}_Z:=\int_D\Delta x\Delta\varphi -\int_D\Psi''_n(y(\omega,t))x\Delta\varphi
  +\int_Dh(\omega,t)\Delta\varphi\,, \qquad 
\]
for $(\omega,t)\in\Omega\times[0,T]$ and $x,\varphi\in Z$. Since $\Psi''_n(y)\in L^\infty(\Omega\times Q)$,
it is not difficult to check that $A_n$ is
progressively measurable, weakly monotone, weakly coercive and linearly bounded.
Moreover, it is clear that the operator
\[
  x\mapsto DB(t,y(\omega,t))x\,, \qquad x\in H\,,
\]
is Lipschitz-continuous and linearly bounded
from $H$ to $\cL^2(U,H)$, uniformly on $\Omega\times[0,T]$
Hence, the approximated problem
admits a unique solution $z_h^n$ such that, setting $\mu_h^n:=-\Delta z_h^n + \Psi_n''(y)z_h^n-h$,
conditions \eqref{lin1}--\eqref{lin4} are satisfied with $\Psi''_n$.

\noindent{\bf Uniform estimates.}
Let us now prove uniform estimates independently of $n$ and pass to the limit as $n\to\infty$.
Noting that $(z_h^n)_D=0$ by (A4), 
It\^o's formula for the square of the $V^*$-norm yields
\begin{align*}
  &\frac12\norm{\nabla\mathcal N(z^n_h)(t)}_H^2 + \int_{Q_t}|\nabla z^n_h|^2
  +\int_{Q_t}\Psi''(y)|z^n_h|^2 = \int_{Q_t}hz_h^n\\
  &\qquad+\frac12\int_0^t\operatorname{Tr}\left(DB(s,y(s))z^n_h(s)^*\circ
  \mathcal N\circ DB(s,y(s))z^n_h(s)\right)\,ds\\
  &\qquad+\int_0^t\left(\mathcal N(z_h^n)(s),DB(s,y(s))z^n_h(s)\right)_H\,dW(s)\,.
\end{align*}
for every $t\in[0,T]$, $\P$-almost surely.
Since $\Psi''\geq-c_1$ implies that $\Psi''_n\geq-c_1$ for every $n\in\enne$,
by the Young inequality we have 
\begin{align*}
  &\frac12\norm{\nabla\mathcal N(z^n_h)(t)}_H^2 + \int_{Q_t}|\nabla z^n_h|^2
  \leq \frac12\int_{Q_t}|h|^2 + \left(\frac12 + c_1\right)\int_{Q_t}|z_h^n|^2\\
  &\qquad+\frac12\int_0^t\operatorname{Tr}\left(DB(s,y(s))z^n_h(s)^*\circ
  \mathcal N\circ DB(s,y(s))z^n_h(s)\right)\,ds\\
  &\qquad+\int_0^t\left(\mathcal N(z_h^n)(s),DB(s,y(s))z^n_h(s)\right)_H\,dW(s)\,.,
\end{align*}
where, by the properties of $\mathcal N$,
\[
  \int_{Q_t}|z_h^n|^2 \leq \delta \int_{Q_t}|\nabla z_h^n|^2 
  +C_\delta\int_0^t\norm{\nabla\mathcal N z_h^n(s)}_H^2\,ds
\]
for every $\delta >0$ and a positive constant $C_\delta >0$.
Hence, choosing $\delta$ sufficiently small, taking 
power $3$, supremum in time and then
expectations, arguing on the right-hand side exactly as in
the proof of uniqueness in Section~\ref{subsec:lin},
we deduce that 
\beq\label{est1_lin}
  \norm{z_h^n}_{L^6(\Omega; C^0([0,T]; V^*)\cap L^2(0,T; V))} 
  \lesssim \norm{h}_{L^6(\Omega; L^2(0,T; H))}
  \qquad\forall\,n\in\enne\,.
\eeq
Now we write It\^o's formula for the square of the $H$-norm, getting
\begin{align*}
  &\frac12\norm{z_h^n(t)}_H^2 + \int_{Q_t}|\Delta z_h^n|^2 =
  \int_{Q_t}\Psi''_n(y)z_h^n\Delta z_h^n - \int_{Q_t}h\Delta z_h^n\\
  &+\int_0^t\norm{DB(s,y(s))z_h^n(s)}^2_{\cL^2(U,H)}\,ds
  +\int_0^t\left(z_h^n(s), DB(s,y(s))z_h^n(s)\right)_H\,dW(s)
\end{align*}
for every $t\in[0,T]$, $\P$-almost surely.
We proceed now similarly to the previous estimate, 
taking supremum in time and expectations.
Using the boundedness of $DB$ together with
the Burkholder-Davis-Gundy and Young inequalities 
on the right-hand side we have in particular that 
\begin{align*}
  &\E\int_0^t\norm{DB(s,y(s))z_h^n(s)}^2_{\cL^2(U,H)}\,ds
  +\E\sup_{r\in[0,t]}\left|\int_0^r\left(z_h^n(s), DB(s,y(s))z_h^n(s)\right)\,dW(s)\right|\\
  &\lesssim \delta\E\norm{z_h^n}^2_{C^0([0,T]; H)} + C_\delta\E\norm{z_h^n}^2_{L^2(0,t; H)}
\end{align*}
for every $\delta >0$ and a certain $C_\delta >0$. Choosing $\delta$
sufficiently small we infer that 
\begin{align*}
  &\norm{z_h^n}_{L^2(\Omega; C^0([0,T]; H))}^2 + 
  \norm{\Delta z_h^n}^2_{L^2(\Omega; L^2(0,T; H))}\\
  &\lesssim
  \norm{h}^2_{L^2(\Omega; L^2(0,T; H))} + \E\int_{Q}|\Psi''_n(y)z_h^n|^2
  +\norm{z_h^n}_{L^2(\Omega; L^2(0,t; H))}^2\,,
\end{align*}
where by (A2), the H\"older inequality 
and the Sobolev embedding $V\embed L^6(D)$ and \eqref{est1_lin},
\begin{align*}
  \E\int_Q|\Psi''_n(y)z_h^n|^2&\lesssim_{c_2} \E\int_Q(1+|y|^4) |z_h^n|^2\\
  &\leq\E \int_0^T\norm{1+|y(s)|^4}_{L^{3/2}(D)}\norm{|z_h^n(s)|^2}_{L^3(D)}\,ds\\
  &\lesssim\E(1+\norm{y}^4_{L^\infty(0,T; V)})\norm{z_h^n}^2_{L^2(0,T; V)}\\
  &\leq(1+\norm{y}^4_{L^6(\Omega; L^\infty(0,T; V)}))\norm{h}^2_{L^6(\Omega; L^2(0,T; H))}\,.
\end{align*}
The estimate \eqref{est1'_state} yields then, thanks to the Gronwall lemma,
\beq\label{est2_lin}
  \norm{z_h^n}_{L^2(\Omega; C^0([0,T]; H)\cap L^2(0,T; Z))} \lesssim 
  \norm{h}_{L^6(\Omega; L^2(0,T; H))}
  \qquad\forall\,n\in\enne\,.
\eeq
By comparison in the equations we also deduce that 
\beq\label{est3_lin}
  \norm{\mu_h^n}_{L^2(\Omega; L^2(0,T; H))}  \lesssim
  \norm{h}_{L^6(\Omega; L^2(0,T; H))}
  \qquad\forall\,n\in\enne\,.
\eeq

\noindent{\bf Passage to the limit.}
By the estimates \eqref{est1_lin}--\eqref{est3_lin}, we deduce that 
there are
\[
  z_h\in L^2(\Omega; C^0([0,T]; H)\cap L^2(0,T; Z))\,,\qquad
  \mu_h \in L^2(\Omega;  L^2(0,T; H))\,,
\]
such that, as $n\to\infty$,
\[
  z_h^{n} \wto z_h \quad\text{in } L^2(\Omega; L^2(0,T; Z))\,,\qquad
  \mu_h^n \wto\mu_h \quad\text{in } L^2(\Omega; L^2(0,T; H))\,.
\]
Moreover, since $DB(\cdot,y)\in\cL(H; \cL^2(U,H))$
and $DB(\cdot,y)^*\in\cL(\cL^2(U,H); H)$, the boundedness of $DB$
ensures that for every $\varphi \in L^2(\Omega; L^2(0,T; \cL^2(U,H)))$
we have that
$DB(\cdot,y)^*\varphi \in L^2(\Omega; L^2(0,T; H))$: hence, the weak
convergence of $(z_h^n)_n$ readily implies that 
\begin{align*}
&\E\int_0^T\left(DB(s,y(s))z_h^n(s),\varphi(s)\right)_{\cL^2(U,H)}\,ds
=\E\int_0^T\left(z_h^n(s),DB(s,y(s))^*\varphi(s)\right)_{H}\,ds\\
&\to\E\int_0^T\left(z_h(s),DB(s,y(s))^*\varphi(s)\right)_{H}\,ds
=\E\int_0^T\left(DB(s,y(s))z_h(s),\varphi(s)\right)_{\cL^2(U,H)}\,ds\,.
\end{align*}
Since $\varphi$ is arbitrary we infer that 
\[
  DB(\cdot, y)z_h^n \wto B(\cdot,y)z_h \quad\text{in } L^2(\Omega; L^2(0,T; \cL^2(U,H)))\,,
\]
hence also, by the linearity and continuity of the stochastic intragral,
\[
  \int_0^\cdot DB(s,y(s))z_h^n(s)\,dW(s) \wto
  \int_0^\cdot DB(s,y(s))z_h^n(s)\,dW(s)
  \quad\text{in } L^2(\Omega; L^2(0,T; H))\,.
\]
It is clear then that these convergences are enough to pass to the limit in 
the variational formulation \eqref{lin3}--\eqref{lin4}, except for the term $\Psi''_n(y)z_h^n$:
let us analyse it explicitly. To this end, note that 
since $\Psi''$ has quadratic growth and $y\in L^6(\Omega; L^\infty(0,T; L^6(D)))$,
we have in particular that $\Psi''(y) \in L^3(\Omega\times(0,T)\times D)$, hence also
\[
  \Psi''_n(y) \to \Psi''(y) \quad\text{in } L^3(\Omega\times(0,T)\times D)\,.
\]
The weak convergence of $(z_h^n)_n$ implies then that 
\[
  \Psi''_n(y)z_h^n \wto \Psi''(y)z_h \quad\text{in } L^{6/5}(\Omega\times(0,T)\times D)\,.
\]
Hence, letting $n\to\infty$ in
\eqref{lin3}--\eqref{lin4} we deduce that $(z_h, \mu_h)$
satisfies the variational formulation of the linearized system.

\subsection{Weak differentiability of the control-to-state map}\label{subsec:diff}
We show here that the map $S$ is weakly G\^ateaux-differentiable in the sense 
specified by Theorem~\ref{thm:3}, and that its weak derivative is the unique 
solution $z_h$ to \eqref{lin1}--\eqref{lin4}.

Let $u, h\in \mathcal U'$ and fix $\eps_0>0$ sufficiently small such that 
$u+\eps h\in\mathcal U'$ for all $\eps\in[-\eps_0,\eps_0]$.
Set also $y:=S(u)$ and $y_h^\eps:=S(u+\eps h)$ 
for any $\eps\in[-\eps0,\eps_0]\setminus\{0\}$, and 
let $z_h$ be the unique solution to the linearized system given by Proposition~\ref{prop:lin}.
Then we have,
in the variational triple $(Z,H,Z^*)$,
\begin{align*}
  &d\left(\frac{y_h^\eps-y}{\eps}\right) - \Delta\left(\frac{w_h^\eps-w}{\eps}\right)\,dt=
  \frac{B(y_h^\eps)-B(y)}\eps\,dW\,,\\
  &\frac{w_h^\eps-w}{\eps} = -\Delta\left(\frac{y_h^\eps-y}{\eps}\right)
  +\frac{\Psi'(y_h^\eps) - \Psi'(y)}\eps - h\,.
\end{align*}
By the continuous dependence properties 
\eqref{dep_cont}--\eqref{dep_cont'}
we have that 
\beq\label{est1_gat}
  \norm{\frac{y_h^\eps-y}{\eps}}_{L^6(\Omega;C^0([0,T]; V^*)\cap L^2(0,T; V))} 
  \lesssim \norm{h}_{L^6(\Omega; L^2(0,T; V^*))}
\eeq
and 
\beq
  \label{est2_gat}
  \norm{\frac{y_h^\eps-y}{\eps}}_{L^2(\Omega; C^0([0,T]; H)\cap L^2(0,T; Z))} 
  \lesssim
  \norm{h}_{L^6(\Omega; L^2(0,T; H))}\,.
\eeq
Furthermore, the mean-value theorem and 
the fact that $\Psi''$ has quadratic growth implies,
by the H\"older inequality and the continuous embedding $V\embed L^6(D)$,
\begin{align*}
  &\E\int_Q\left|\frac{\Psi'(y_h^\eps) - \Psi'(y)}\eps\right|^2\leq
  \E\int_Q\int_0^1|\Psi''(y + \sigma(y_h^\eps-y))|^2\left|\frac{y_h^\eps-y}\eps\right|^2\,d\sigma\\
  &\lesssim \E\int_Q\left(1+|y|^4 + |y_h^\eps|^4\right)\left|\frac{y_h^\eps-y}\eps\right|^2\\
  &\lesssim\left(1+ \norm{y}^4_{L^6(\Omega; L^\infty(0,T; V))}
  + \norm{y_h^\eps}^4_{L^6(\Omega; L^\infty(0,T; V))}\right)
  \norm{\frac{y_h^\eps-y}\eps}_{L^6(\Omega; L^2(0,T; V))}^2\,,
\end{align*}
so that \eqref{est1'_state} and \eqref{est1_gat} imply that 
\beq
  \label{est3_gat}
  \norm{\frac{\Psi'(y_h^\eps) - \Psi'(y)}\eps}_{L^2(\Omega; L^2(0,T; H))}\lesssim
  \norm{h}_{L^6(\Omega; L^2(0,T; H))}\,.
\eeq
Moreover, the Lipschitz-continuity of $B$ and \eqref{est2_gat} ensures that 
\[
  \norm{\frac{B(\cdot, y_h^\eps)-B(\cdot, y)}\eps}_{L^2(\Omega; C^0([0,T]; \cL^2(U,H)))}
  \lesssim
  \norm{h}_{L^6(\Omega; L^2(0,T; H))}\,,
\]
from which 
\beq\label{est4_gat}
\norm{\int_0^\cdot\frac{B(s,y_h^\eps(s))-B(s,y(s))}\eps\,dW(s)}_{L^2(\Omega; C^0([0,T]; H))}
\lesssim\norm{h}_{L^6(\Omega; L^2(0,T; H))}\,.
\eeq
Hence by comparison in the equation we also have that
\beq\label{est5_gat}
  \norm{\frac{w_h^\eps-w}{\eps}}_{L^2(\Omega; L^2(0,T; H))}
 \lesssim
  \norm{h}_{L^6(\Omega; L^2(0,T; V))}\,.
\eeq
Let us pass to the limit as $\eps\searrow0$. From the 
estimates \eqref{est1_gat}--\eqref{est5_gat}
we deduce that there are
\[
  z_h \in L^2\left(\Omega; L^\infty(0,T; H)\cap L^2(0,T; Z)\right)\,, \quad
  \mu_h \in L^2_{\mathcal P}(\Omega; L^2(0,T; H))
\]
such that, as $\eps\searrow0$,
\begin{align}
  \label{conv1_incr}
  \frac{y_h^\eps-y}{\eps} \wto z_h \qquad&\text{in } 
  L^2\left(\Omega; L^p(0,T; H)\cap L^2(0,T; Z)\right) \quad\forall\,p\in[1,+\infty)\,,\\
  \label{conv2_incr}
  \frac{w_h^\eps-w}{\eps} \wto \mu_h \qquad&\text{in } L^2\left(\Omega; L^2(0,T; H)\right)\,.
\end{align}

Moreover, note that 
\begin{align*}
  &\frac{\Psi'(y_h^\eps) - \Psi'(y)}\eps - \Psi''(y)z_h\\
  &=\frac{\Psi'(y_h^\eps) - \Psi'(y) - \Psi''(y)(y_h^\eps-y)}\eps
  +\Psi''(y)\left(\frac{y_h^\eps - y}{\eps} - z_h\right)\\
  &=\int_0^1\left(\Psi''(y + r(y_h^\eps-y)) - \Psi''(y)\right)\frac{y_h^\eps - y}{\eps}\,dr
  +\Psi''(y)\left(\frac{y_h^\eps - y}{\eps} - z_h\right)\,.
\end{align*}
Since $\Psi''(y)\in L^3(\Omega\times(0,T)\times D)$,
the weak convergences proved above imply that
\[
  \Psi''(y)\left(\frac{y_h^\eps - y}{\eps} - z_h\right) \wto 0
  \quad\text{in } L^{6/5}(\Omega\times(0,T)\times D)\,.
\]
Let us show that also the first term goes to $0$.
To this end, note that since
\[
  \norm{y_h^\eps-y}_{L^2(\Omega; C^0([0,T]; H)\cap L^2(0,T; Z))}\lesssim \eps
  \norm{h}_{L^6(\Omega; L^2(0,T; H))}\to 0\,,
\]
by continuity of $\Psi''$, we have, along a subsequence,
\[
  \Psi''(y + r(y_h^\eps-y)) - \Psi''(y) \to 0 \quad\forall\,r\in[0,1]\,,\quad
  \text{a.e.~in } \Omega\times(0,T)\times D\,.
\]
Moreover, 
since $\Psi''$ has quadratic growth, we deduce that
\[
  \left|\int_0^1(\Psi''(y + r(y_h^\eps-y)) - \Psi''(y))\,dr\right|\lesssim
  1+|y_h^\eps|^2 + |y|^2\,,
\]
where the right hand side is bounded in $L^3(\Omega\times(0,T)\times D)$
because $y$ and $(y_h^\eps)_\eps$ are bounded in $L^6(\Omega; L^\infty(0,T; V))$ 
by Theorem~\ref{thm:1} and $V\embed L^6(D)$.
Consequently, 
\[
  \int_0^1\left(\Psi''(y + r(y_h^\eps-y)) - \Psi''(y)\right)\,dr \to 0 \quad\text{in }
  L^p(\Omega\times(0,T)\times D) \quad\forall\,p\in[2,3)\,.
\]
In particular, we deduce that 
\[
  \int_0^1\left(\Psi''(y + r(y_h^\eps-y)) - \Psi''(y)\right)\frac{y_h^\eps - y}{\eps}\,dr
  \wto 0 \quad\text{in } L^{p}(\Omega\times(0,T)\times D)
\]
for all $p\in[1,6/5)$, from which 
\beq\label{conv3_incr}
  \frac{\Psi'(y_h^\eps) - \Psi'(y)}\eps \wto 
  \Psi''(y)z_h \quad\text{in } L^{p}(\Omega\times(0,T)\times D)
  \quad\forall\,p\in[1,6/5)\,.
\eeq

Let us show the convergence of the stochastic integrals. To this end, note that 
\begin{align*}
  &\frac{B(\cdot, y^\eps_h)-B(\cdot, y)}{\eps} - DB(\cdot, y)z_h\\
  &=\frac{B(\cdot, y^\eps_h)-B(\cdot, y) - DB(\cdot, y)(y^\eps-y)}{\eps} 
  +DB(\cdot, y)\left(\frac{y^\eps_h-y}{\eps} - z_h\right)\\
  &=\int_0^1\left(DB(\cdot, y+r(y_h^\eps-y))-DB(\cdot, y)\right)\frac{y_h^\eps-y}\eps\,dr
  +DB(\cdot, y)\left(\frac{y^\eps_h-y}{\eps} - z_h\right)\,.
\end{align*}
The weak convergences proved above, the linearity of $DB(\cdot, y)$, 
the boundedness of $DB$ and the dominated convergence theorem yields
\[
  DB(\cdot, y)\left(\frac{y^\eps_h-y}{\eps} - z_h\right)\wto 0
  \quad\text{in } L^2(\Omega; L^2(0,T; \cL^2(U,H)))\,.
\]
Moreover, since 
$DB(t,\cdot)\in C^0(H; \cL(H,\cL^2(U,H)))$ by assumption (A5),
recalling also that $y_h^\eps\to y$ in $L^2(\Omega; C^0([0,T]; H))$,
by the dominated convergence theorem we have that 
\[
  \int_0^1\left(DB(\cdot, y+r(y_h^\eps-y))-DB(\cdot, y)\right)\,dr \to 0 \quad\text{in }
  L^p(\Omega; L^p(0,T; \cL(H;\cL^2(U;H))))
\]
for every $p\in[2,\infty)$.
Since $\frac{y_h^\eps-y}{\eps}\wto z_h$ in $L^2(\Omega; L^p(0,T; H))$, we deduce
in particular that 
\[
  \int_0^1\left(DB(\cdot, y+r(y_h^\eps-y))-DB(\cdot, y)\right)\frac{y_h^\eps-y}\eps\,dr
  \wto 0 \quad\text{in } L^p(\Omega; L^2(0,T; \cL^2(U,H)))\,,
\] 
Consequently, taking this information into account, we have
\[
  \frac{B(\cdot, y^\eps_h)-B(\cdot, y)}{\eps} \wto DB(\cdot, y)z_h
  \quad\text{in } L^p(\Omega; L^2(0,T;\cL^2(U;H)))\quad
  \forall\,p\in[1,2)\,,
\]
from which
\beq\label{conv5_incr}
  \int_0^\cdot\frac{B(s, y^\eps_h(s))-B(s, y(s))}{\eps}\,dW(s) 
  \wto \int_0^\cdot DB(s, y(s))z_h(s)\,dW(s)
\eeq
in $L^p(\Omega; L^2(0,T;H))$ for all $p\in[1,2)$.

Hence, letting $\eps\to0$ in the variational formulation we deduce that $(z_h, \mu_h)$
satisfy the linearized system \eqref{lin1}--\eqref{lin4}.
Since we have already proved uniqueness for such system in the previous section, 
we deduce that $(z_h, \mu_h)$ is the unique solution to \eqref{lin1}--\eqref{lin4}.

\subsection{First-order necessary conditions for optimality}
We prove here the version of first-order necessary optimality conditions
contained in Theorem~\ref{thm:4}.

Let $\bar u\in \mathcal U$ be an optimal control and set $\bar u:= S(\bar u)$.
For every $v\in \mathcal U$ let us define $h:=v-u$,
and $y^\eps_h:=S(\bar u + \eps h)$ for every $\eps>0$.
Since $\mathcal U$ is convex, we have that $u+\eps(v-u)\in \mathcal U$
for all $\eps\in[0,1]$: hence,
by definition of optimal control we have that $\tilde J(\bar u)\leq\tilde J(\bar u+\eps h)$, 
which may be rewritten
\[
J(\bar y, \bar u) \leq \frac{\alpha_1}2\E\int_Q|y_h^\eps - x_Q|^2 
+ \frac{\alpha_2}2\int_D|y_h^\eps(T)-x_T|^2 +\frac{\alpha_3}{2}\E\int_Q|\bar u+\eps h|^2\,.
\]
Using the definition of $J(\bar y, \bar u)$ and rearranging the terms we have 
\begin{align*}
0&\leq \frac{\alpha_1}2\E\int_Q\left(|y_h^\eps|^2 - |\bar y|^2 - 2(y_h^\eps-\bar y)x_Q\right)\\
&+\frac{\alpha_2}2\E\int_D\left(|y_h^\eps(T)|^2 - |\bar y(T)|^2 - 2(y_h^\eps-\bar y)(T)x_T\right)
+\frac{\alpha_3}2\E\int_Q\left(\eps^2|h|^2 + 2\eps \bar u h\right)\,.
\end{align*}
Since the functions $x\mapsto \E\int_Q|x|^2$ and $x\mapsto \E\int_D|x|^2$
are Fr\'echet-differentiable in $L^2(\Omega\times Q)$ and $L^2(\Omega\times D)$, respectively, 
dividing by $\eps$ we get 
\begin{align*}
  0&\leq \alpha_1\E\int_Q\left(\int_0^1(\bar y + \sigma(y_h^\eps-\bar y))\,d\sigma-x_Q\right)
  \frac{y_h^\eps-\bar y}{\eps}\\
  &+\alpha_2\E\int_D\left(\int_0^1(\bar y + \sigma(y_h^\eps-\bar y))(T)\,d\sigma-x_T\right)
  \frac{y_h^\eps-\bar y}{\eps}(T)\\
  &+\alpha_3\E\int_Q\bar u h + \frac{\alpha_3}2\eps\norm{h}^2_{L^2(\Omega\times Q)}\,.
\end{align*}
Since $\bar u+\eps h \to \bar u$ in $L^6(\Omega; L^2(0,T; V))$ as $\eps\searrow0$, 
we deduce from \eqref{dep_cont}--\eqref{dep_cont'},
the definition of $\mathcal U$
and the dominated convergence theorem that 
\begin{align*}
  \int_0^1(\bar y + \sigma(y_h^\eps-\bar y))\,d\sigma-x_Q \to \bar y - x_Q \qquad&\text{in }
  L^6(\Omega; L^2(0,T; H))\,,\\
  \int_0^1(\bar y + \sigma(y_h^\eps-\bar y))(T)\,d\sigma-x_T \to \bar y(T) - x_T \qquad&\text{in }
  L^2(\Omega; H)\,.
\end{align*}
Furthermore, by Theorem~\ref{thm:3} we know that 
\begin{align*}
  \frac{y_h^\eps-\bar y}{\eps} \to z_h \qquad&\text{in }
  L^{6/5}(\Omega; L^2(0,T; H))\,,\\
  \frac{y_h^\eps-\bar y}{\eps}(T) \wto z_h(T) \qquad&\text{in }
  L^2(\Omega; H)\,,
\end{align*}
so that letting $\eps\searrow0$ in the last inequality Theorem~\ref{thm:4} is proved.

\section{The adjoint problem}
\label{sec:ad}
\setcounter{equation}{0}

In this section we study the adjoint problem \eqref{eq:1_ad}--\eqref{eq:fin_ad}
in terms of existence and uniqueness of solutions. Moreover, we 
prove the refined version of first-order necessary optimality conditions
contained in Theorem~\ref{thm:5}.

\subsection{Existence-uniqueness of the adjoint problem}
\label{subsec:ad}
We prove here Proposition~\ref{prop:ad}. 
Let $u\in \mathcal U'$ and
$y:=S(u)$.

\noindent{\bf Uniqueness.}
First of all we prove uniqueness of solutions. Let $(p_i, \tilde p_i, q_i)$
satisfy \eqref{ad1}--\eqref{ad4} for $i=1,2$: taking the difference
of the respective equations we have,
setting $p:=p_1-p_2$, $\tilde p:=\tilde p_1-\tilde p_2$, and $q:=q_1-q_2$,
\[
  -dp - \Delta\tilde p\,dt +\Psi''(y)\tilde p\,dt=
  DB(y)^*q\,dt
  -q\,dW\,,\qquad \tilde p=-\Delta p\,.
\]
It\^o's formula for $\frac12\norm{\nabla p}_H^2$ yields then
\begin{align*}
  &\frac12\E\norm{\nabla p(t)}_H^2
  +\E\int_t^T\!\!\int_D|\nabla\tilde p(s)|^2\,ds
  +\E\int_t^T\!\!\int_D\Psi''(y(s))|\tilde p(s)|^2\\
  &\qquad+\frac12\E\int_t^T\norm{\nabla q(s)}^2_{\cL^2(U,H)}\,ds
  =\E\int_t^T\left(q(s), DB(s,y(s))\tilde p(s)\right)_{\cL^2(U,H)}\,ds\,.
\end{align*}
Recalling assumption (A4), we have that $DB(\cdot,y)\tilde p\in \cL^2(U,H_0)$, 
so that 
\[
  \left(q, DB(\cdot,y)\tilde p\right)_{\cL^2(U,H)} = \left(q-q_D, DB(\cdot,y)\tilde p\right)_{\cL^2(U,H)}\,.
\]
Taking into account (A2) and noting that $\tilde p_D=0$,
we get,
by the Young and Poincar\'e inequalities and \eqref{comp_ineq},
\begin{align*}
  &\frac12\E\norm{\nabla p(t)}_H^2
  +\E\int_t^T\!\!\int_D|\nabla\tilde p(s)|^2\,ds
  +\frac12\E\int_t^T\norm{\nabla q(s)}^2_{\cL^2(U,H)}\,ds \\
  &\leq c_1\E\int_t^T\!\!\int_D|\tilde p(s)|^2\,ds
  +L_B\E\int_t^T\norm{(q-q_D)(s)}_{\cL^2(U,H)}(1+\norm{\tilde p(s)}_H)\,ds\\
  &\leq\sigma\E\int_t^T\!\!\int_D|\nabla \tilde p(s)|^2\,ds
  +C_\sigma\E\int_t^T\norm{\nabla\mathcal N \tilde p(s)}_H^2\,ds
  +\sigma\E\int_t^T\norm{\nabla q(s)}_{\cL^2(U,H)}^2\,ds\\
  &\lesssim\sigma\E\int_t^T\!\!\int_D|\nabla \tilde p(s)|^2\,ds
  +\sigma\E\int_t^T\norm{\nabla q(s)}_{\cL^2(U,H)}
  +C_\sigma\E\int_t^T\norm{\nabla p(s)}_H^2\,ds
\end{align*}
for every $\sigma>0$ for a certain $C_\sigma>0$.
Choosing $\sigma$ sufficiently small and applying the Gronwall lemma yields then
$\nabla \tilde p=0$, from which 
$\tilde p=0$ since $\tilde p_D=0$. Since $\tilde p=-\Delta p$,
we infer that $-\Delta p=0$, from which $p_1-(p_1)_D=p_2-(p_2)_D$.

\noindent{\bf Approximation.}
Let us prove now existence of solution to the BSPDE \eqref{eq:1_ad}--\eqref{eq:fin_ad}.
We perform the same approximation that we used for the 
linearized system in Section~\ref{subsec:lin}, and we consider  for every $n\in\enne$
the approximated problem 
\begin{align*}
    \tilde p_n = -\Delta p_n\qquad&\text{in } Q\,,\\
    -dp_n - \Delta \tilde p_n\,dt +\Psi''_n(y)\tilde p_n\,dt
    = \alpha_1(y-x_Q)\,dt 
    +DB(y)^*q_n\,dt-q_n\,dW \qquad&\text{in } Q\,,\\
    \partial_{\bf n} p_n = \partial_{\bf n} \tilde p_n = 0 \qquad&\text{in } \Sigma\,,\\
     p_n(T)=\alpha_2(y(T)-x_T) \qquad&\text{in } D\,,
\end{align*}
where $\Psi''_n:=T_n\circ\Psi''$ and $T_n:\erre\to\erre$ is the truncation operator at level $n$.
The variational formulation of the approximated problem is given by 
\begin{align*}
  \int_Dp_n(t)\varphi &+ \int_t^T\!\!\int_D\Delta p_n(s)\Delta\varphi\,ds
  -\int_t^T\!\!\int_D\Psi_n''(y(s))\Delta p_n(s)\varphi\,ds\\
  &=\alpha_2\int_D(y(T)-x_T)\varphi + \alpha_1\int_t^T\!\!\int_D(y-x_Q)(s)\varphi\,ds\\
  &+\int_t^T\left(DB(s,y(s))^*q_n(s),\varphi\right)_H\,ds
  -\int_D\left(\int_t^Tq_n(s)\,dW(s)\right)\varphi
\end{align*}
for every $\varphi\in Z$, $\P$-almost surely, for every $t\in[0,T]$. Hence, we introduce the operator
$A_n^*:\Omega\times[0,T]\times Z \to Z^*$ as
\[
  \ip{A_n^*(\omega,t,x)}{\varphi}_{Z}:=\int_D\Delta x\Delta\varphi
  -\int_D\Psi_n''(y(\omega,t))\Delta x\varphi
\]
and note that since $\Psi''_n(y)\in L^\infty(\Omega\times Q)$, then $A_n^*$
is progressively measurable, weakly monotone, weakly coercive
and linearly bounded. Moreover, 
the operator $DB(\cdot,y)^*$ is uniformly bounded
in $\Omega\times[0,T]$ be (A4). Hence, by the classical variational approach to BSPDEs
(see \cite[Sec.~3]{du-meng2}) the approximated problem admits 
a unique solution $(p_n, \tilde p_n, q_n)$ with 
\begin{gather*}
  p_n \in L^2_{\mathcal P}(\Omega; C^0([0,T]; H))\cap L^2_{\mathcal P}(\Omega; L^2(0,T; Z))\,, \\
  \tilde p_n \in L^2_{\mathcal P}(\Omega; C^0([0,T]; Z^*))\cap L^2_{\mathcal P}(\Omega; L^2(0,T; H))\,,\\
  q_n \in L^2_\mathcal P(\Omega; L^2(0,T; \cL^2(U,H)))\,.
\end{gather*}
Moreover, by assumption we have $\alpha_ 2 x_T \in L^2(\Omega,\cF_T,\P; V)$, 
while by Theorem~\ref{thm:1} we know that $y \in L^2(\Omega; C^0([0,T]; H)\cap L^\infty(0,T; V))$,
so that $y$ is weakly continuous in $V$ and $y(T) \in L^2(\Omega,\cF_T,\P; V)$.
Consequently, we have that $\alpha_2(y(T)-x_T) \in L^2(\Omega,\cF_T,\P; V)$, and this 
ensures a further regularity on $(p_n, \tilde p_n, q_n)$, namely
\begin{gather*}
  p_n \in L^2_{\mathcal P}(\Omega; C^0([0,T]; V))\cap 
  L^2_{\mathcal P}(\Omega; L^2(0,T; H^3(D)))\,, \\
  \tilde p_n \in L^2_{\mathcal P}(\Omega; C^0([0,T]; V^*))\cap L^2_{\mathcal P}(\Omega; L^2(0,T; V))\,,\\
  q_n \in L^2_\mathcal P(\Omega; L^2(0,T; \cL^2(U,V)))\,.
\end{gather*}
In order to prove this, one should perform a further approximation on the problem
depending on a further parameter (let us say $k$, for example), 
write It\^o's formula for $\frac12\norm{\nabla p^k_n}_H^2$ and then pass to the limit
as $k\to \infty$. Since the procedure is quite standard,
to avoid heavy notations we shall proceed formally
writing It\^o's formula for $\frac12\norm{\nabla p_n}_H^2$: this reads
\beq\label{ito_pn}
\begin{split}
  &\frac12\norm{\nabla p_n(t)}_H^2
  +\int_t^T\!\!\int_D|\nabla\Delta p_n(s)|^2\,ds +
  \int_t^T\!\!\int_D\Psi''_n(y(s))|\Delta p_n(s)|^2\,ds\\
  &\qquad+\frac12\int_t^T\norm{\nabla q_n(s)}^2_{\cL^2(U,H)}\,ds
  -\int_t^T\left(\Delta p_n(s), q_n(s)\right)_H\,dW(s)\\
  &=\frac{\alpha_2^2}2\norm{\nabla(y(T)-x_T)}_H^2
  -\alpha_1\int_t^T\!\!\int_D(y-x_Q)(s)\Delta p_n(s)\,ds\\
  &\qquad+\int_t^T\left(DB(s,y(s))^*q_n(s),\tilde p_n(s)\right)_H\,ds\,.
\end{split}
\eeq
Since $\Psi''_n\in L^\infty(\Omega\times Q)$ (recall that here $n$ is fixed)
and we already know that $p_n\in L^2(\Omega; L^2(0,T; Z))$, the desired regularity
follows by a classical procedure based on the Burkholder-Davis-Gundy inequality.
The regularity for $\tilde p_n$ follows then by comparison.

\noindent{\bf First estimate.}
We now prove uniform estimates independently of $n$
and pass to the limit as $n\to\infty$.
First of all, taking expectations in \eqref{ito_pn},
and performing the same computations as in 
the proof of uniqueness at the beginning of Section~\ref{subsec:ad} yields
\begin{align*}
  &\E\norm{\tilde p_n(t)}_{V^*}^2
  +\E\int_t^T\!\!\int_D|\nabla\tilde p_n(s)|^2\,ds
  +\E\int_t^T\norm{\nabla q_n(s)}^2_{\cL^2(U,H)}\,ds\\
  &\lesssim_{c_1,L_B}
  \norm{y(T)}^2_{L^2(\Omega; V)} + \norm{\alpha_2 x_T}^2_{L^2(\Omega; V)}
  +\alpha_1^2\E\int_Q|y-x_Q|^2\\
  &\qquad+\E\int_t^T\!\!\int_D|\tilde p_n(s)|^2\,ds
  +\E\int_t^T\norm{\nabla q_n(s)}_{\cL^2(U,H)}^2\,ds\,.
\end{align*}
Recalling that $\norm{\nabla p_n}_H=\norm{\nabla \mathcal N\tilde p_n}_H
\lesssim\norm{\tilde p_n}_{V^*}$, using the compactness inequality \eqref{comp_ineq}
on the right-hand side yields, by the Gronwall lemma,
\[
  \norm{\nabla q_n}_{L^2(\Omega; L^2(0,T; \cL^2(U,H)))}\leq c\,.
\]
Hence, going back again in It\^o's formula \eqref{ito_pn}, we now take supremum in time
and then expectations: we estimate the stochastic integral using the Burkholder-Davis Gundy inequality
and integration by parts as (see.~e.g.~\cite[Lem.~4.3]{mar-scar-diss})
\begin{align*}
  &\E\sup_{t\in[0,T]}\left|\int_t^T\left(\Delta p_n(s), q_n(s)\right)_H\,dW(s)\right|\\
  &\qquad\lesssim\eps \E\norm{\nabla p_n}_{C^0([0,T]; H)}^2 
  +C_\eps\E\norm{\nabla q_n}^2_{L^2(0,T; \cL^2(U,H))}
\end{align*}
for every $\eps>0$, so that choosing $\eps$ sufficiently small, rearranging the terms 
and recalling the estimate just proved on $(\nabla q_n)_n$ yields,
for a positive constant $c$ independent of $n$,
\beq\label{est1_ad}
  \norm{\tilde p_n}_{L^2(\Omega; C^0([0,T]; V^*)))\cap L^2(\Omega; L^2(0,T; V))} 
  +\norm{\nabla q_n}_{L^2(\Omega; L^2(0,T; \cL^2(U,H)))}\leq c\,.
\eeq

\noindent{\bf Second estimate.} 
We write It\^o's formula for $(\frac12\norm{\nabla p_n}^2_H)^3$, getting
\beq\label{ito_pn6}
\begin{split}
  &\frac18\norm{\nabla p_n(t)}_H^6 + 
  \frac34\int_t^T\norm{\nabla p_n(s)}_H^4\norm{\nabla \tilde p_n(s)}_H^2\,ds\\
  &+\frac34\int_t^T\norm{\nabla p_n(s)}_H^4\!\!\int_D\Psi''(y)|\tilde p_n(s)|^2\,ds
  +\frac38\int_t^T\norm{\nabla p_n(s)}_H^4\norm{\nabla q_n(s)}_{\cL^2(U,H)}^2\,ds\\
  &+\frac32\int_t^T\norm{\nabla p_n(s)}_H^2
  \norm{\ip{-\Delta p_n(s)}{q_n(s)}_V}^2_{\cL^2(U,\erre)}\,ds\\
  &-\frac34\int_t^T\norm{\nabla p_n(s)}_H^4\left(\Delta p_n(s), q_n(s)\right)_H\,dW(s)\\
  &=\frac18\norm{\alpha_2\nabla(y(T)-x_T)}_H^6
  +\frac34\alpha_1\int_t^T\!\!\int_D\norm{\nabla p_n(s)}_H^4(y-x_Q)(s)\tilde p_n(s)\,ds\\
  &+\frac34\int_t^T\norm{\nabla p_n(s)}_H^4
  \left(q_n(s), DB(s,y(s))\tilde p_n(s)\right)_{\cL^2(U,H)}\,ds\,,
\end{split}
\eeq
By the H\"older and Young inequalities and the definition of $\tilde p_n$,
for all $\delta>0$ we have
\begin{align*}
&\alpha_1\int_t^T\!\!\int_D\norm{\nabla p_n(s)}_H^4(y-x_Q)(s)\tilde p_n(s)\,ds\\
&\leq\alpha_1\int_t^T\norm{\nabla p_n(s)}_H^4\norm{\tilde p_n(s)}_H\norm{(y-x_Q)(s)}_H\,ds\\
&\leq \delta\int_t^T\norm{\nabla p_n(s)}_H^4\norm{\nabla \tilde p_n(s)}_H^2\,ds+
\norm{\alpha_1(y-x_Q)}^6_{L^6(0,T; H)}+
C_\delta\int_t^T\norm{\nabla p_n(s)}_H^6\,ds\,.
\end{align*}
Moreover, 
recalling that $\Psi''\geq-c_1$ ad $\norm{\tilde p_n}_{V^*}\lesssim\norm{\nabla p_n}_H$, 
for every $\delta>0$ we have
\begin{align*}
  &-\E\int_t^T\norm{\nabla p_n(s)}_H^4\!\!\int_D\Psi''(y(s))|\tilde p_n(s)|^2\,ds
  \leq c_1\E\int_t^T\norm{\nabla p_n(s)}_H^4\norm{\tilde p_n(s)}^2_H\,ds\\
  &\leq \delta \E\int_t^T\norm{\nabla p_n(s)}_H^4\norm{\nabla \tilde p_n(s)}^2_H\,ds
  +C_\delta\int_t^T\E\norm{\nabla p_n(s)}_H^6\,ds\,.
\end{align*}
Similarly, by assumptions (A4)--(A5), recalling that 
$(DB(\cdot, y)\tilde p_n)_D=0$ and writing $q=q-q_D+q_D$, and
arguing as in the the proof of \eqref{est1_ad} we have
\begin{align*}
  &\E\int_t^T\norm{\nabla p_n(s)}_H^4
  \left(q_n(s), DB(s,y(s))\tilde p_n(s)\right)_{\cL^2(U,H)}\,ds\\
  &\lesssim_{L_B}1+
  \delta \E\int_t^T\norm{\nabla p_n(s)}_H^4\norm{\nabla q_n(s)}_{\cL^2(U,H)}^2\,ds\\
  &\qquad
  +\delta \E\int_t^T\norm{\nabla p_n(s)}_H^4\norm{\nabla\tilde p_n(s)}_H^2\,ds
  +C_\delta\int_t^T\norm{\nabla p_n(s)}_H^6\,ds\,.
\end{align*}
Hence, taking expectations in \eqref{ito_pn6}, 
recalling the assumptions on $x_T$ and $x_Q$ and that $y\in L^6(\Omega; L^\infty(0,T; V))$, 
choosing $\delta>0$ sufficiently small,
the Gronwall lemma yields
\[
  \norm{\nabla p_n}^6_{C^0([0,T]; L^6(\Omega; H)))} + 
  \E\int_0^T\norm{\nabla p_n(s)}_H^4\norm{\nabla q_n(s)}_{\cL^2(U,H)}^2\,ds \leq c\,.
\]
At this point, we go back to \eqref{ito_pn6}, take supremum in time and then 
expectations: estimating the stochastic integral through the Burkolder-Davis-Gundy inequality as
\begin{align*}
  &\E\sup_{t\in[0,T]}\left|\int_t^T\norm{\nabla p_n(s)}_H^4\left(\Delta p_n(s), q_n(s)\right)_H\,dW(s)\right|\\
  &\qquad
  \lesssim\E\left(\int_0^T\norm{\nabla p_n(s)}_H^{10}\norm{\nabla q_n(s)}^2_{\cL^2(U,H)}\,ds\right)^{1/2}\\
  &\qquad\leq\E\norm{\nabla p_n}^3_{C^0([0,T]; H)}
  \left(\int_0^T\norm{\nabla p_n(s)}^4_{H}\norm{\nabla q_n(s)}^2_{\cL^2(U,H)}\,ds\right)^{1/2}\\
  &\qquad\leq\delta\E\norm{\nabla p_n}^6_{C^0([0,T]; H)} + C_\delta
  \E\int_0^T\norm{\nabla p_n(s)}_H^4\norm{\nabla q_n(s)}_{\cL^2(U,H)}^2\,ds\,.
\end{align*}
Choosing $\delta>0$, rearranging the terms and taking into account the estimate already proved, we get
\beq
  \label{est1'_ad}
  \norm{\nabla p_n}^6_{L^6(\Omega; C^0([0,T]; H))}
  +\E\int_0^T\norm{\nabla p_n(s)}_H^4\norm{\nabla q_n(s)}_{\cL^2(U,H)}^2\,ds \leq c\,.
\eeq
Finally, we go back to \eqref{ito_pn}, take the $3$rd-power and then expectations:
using again (A2) and the Young inequality we get
\begin{align*}
  &\E\sup_{r\in[t,T]}\norm{\nabla p_n(r)}_H^6
  +\E\norm{\nabla\tilde p_n}_{L^2(t,T; H)}^6 +
  \E\norm{\nabla q_n}^6_{L^2(t;T;\cL^2(U,H))}\\
  &\lesssim_{c_1,L_B} \E\norm{\alpha_2\nabla (y(T)-x_T)}^6_H
  +\E\norm{\tilde p_n}_{L^2(t,T;H)}^6
  +\E\norm{\alpha_1(y-x_Q)}^6_{L^2(0,T; H)}\\
  &\qquad+\E\sup_{r\in[0,T]}\left|\int_r^T\left(\Delta p_n(s), q_n(s)\right)_H\,dW(s)\right|^3\,,
\end{align*}
where the last term is estimated thanks to the Burkholder-Davis-Gundy inequality by
\begin{align*}
  &\E\left(\int_t^T\norm{\nabla p_n(s)}_H^2\norm{\nabla q_n(s)}_{\cL^2(U,H)}^2\,ds\right)^{3/2}\\
  &\qquad\leq\sigma\E\norm{\nabla q_n}_{L^2(t,T; \cL^2(U,H))}^6
  +C_\sigma\E\norm{\nabla p_n}_{C^0([0,T]; H)}^6
\end{align*}
for every $\sigma>0$ and for a certain $C_\sigma>0$. Hence, noting also that 
\[
  \E\norm{\tilde p_n}_{L^2(t,T;H)}^6\leq\sigma\E\norm{\nabla \tilde p_n}^6_{L^2(t,T;H)}
  +\tilde C_\sigma\E\norm{\nabla p_n}_{L^2(t,T; H)}^6
\]
for a certain $\tilde C_\sigma>0$, 
choosing $\sigma>0$ sufficiently small, rearranging the terms and taking \eqref{est1'_ad} into account,
by the Gronwall lemma we deduce that 
\beq
  \label{est1''_ad}
  \norm{\tilde p_n}_{L^6(\Omega; L^2(0,T; V))} + 
  \norm{\nabla q_n}_{L^6(\Omega; L^2(0,T; \cL^2(U,H)))}\leq c\,.
\eeq

\noindent{\bf Third estimate.}
We write It\^o's formula for $\frac12\norm{p_n}_H^2$, getting
for every $t\in[0,T]$, $\P$-a.s.
\begin{align*}
  &\frac12\norm{p_n(t)}_H^2
  +\int_t^T\!\!\int_D|\nabla p_n(s)|^2\,ds
  +\int_t^T\!\!\int_D\Psi''_n(y(s))\tilde p_n(s) p_n(s)\,ds\\
  &\qquad+\frac12\int_t^T\norm{q_n(s)}^2_{\cL^2(U,H)}\,ds
  +\int_t^T\left(p_n(s), q_n(s)\right)_H\,dW(s)\\
  &=\frac{\alpha_2^2}2\norm{y(T)-x_T}_H^2
  +\alpha_1\int_t^T\!\!\int_D(y-x_Q)(s) p_n(s)\,ds\\
  &\qquad+\int_t^T\left(q_n(s), DB(s,y(s))p_n(s)\right)_{\cL^2(U,H)}\,ds\,.
\end{align*}
Taking expectations,
using the Young inequality, the boundedness of $DB$ and the 
estimates \eqref{est1_ad}--\eqref{est1''_ad}, we infer that, for every $t\in[0,T]$,
\[
  \E\norm{p_n(t)}_H^2
  +\E\int_t^T\!\!\int_D|\nabla p_n(s)|^2\,ds
  +\E\int_t^T\norm{q_n(s)}^2_{\cL^2(U,H)}\,ds
  \lesssim 1
  +\E\int_Q|\Psi''_n(y)\tilde p_n|^2\,,
\]
where the implicit constant is independent of $n$.
Now, by (A2) and the H\"older and Young inequalities, we have
\begin{align*}
  \E\int_Q|\Psi''(y)\tilde p_n|^2
  &\lesssim\E\int_Q|\tilde p_n|^2 + \E\int_Q|y|^4|\tilde p_n|^2\\
  &\lesssim
  \E\int_Q|\tilde p_n|^2 + \E\norm{y}_{L^\infty(0,T;V)}^4\norm{\tilde p_n}_{L^2(0,T; V)}^2\\
  &\leq \E\norm{\tilde p_n}^2_{L^2(0,T; H)}
  +\E\norm{y}_{L^\infty(0,T; V)}^6 + \E\norm{\tilde p_n}_{L^2(0,T; V)}^6\,,
\end{align*}
so that by \eqref{est1_ad} and \eqref{est1''_ad} we get 
\beq\label{est2_ad}
  \norm{q_n}_{L^{2}(\Omega; L^2(0,T; \cL^2(U,H)))}+
   \norm{\Psi''_n(y)\tilde p_n}_{L^{2}(\Omega; L^2(0,T; H))}\leq c\,.
\eeq
Now we go back to It\^o's formula for $\frac12\norm{p_n}_H^2$: instead of taking expectations
straight away, 
we take at first supremum in time and then expectations, getting
Performing the usual computation as before using the Young inequality we get, for every $t\in[0,T]$,
\[
  \E\sup_{r\in[t,T]}\norm{p_n(t)}_H^2
  \lesssim 1+  \E\sup_{r\in[t,T]}\left|\int_r^T\left(p_n(s), q_n(s)\right)_H\,dW(s)\right|\,.
\]
The Burkholder-Davis-Gundy and Young inequalities ensure that, 
for every $\delta >0$, 
\[
\E\sup_{r\in[t,T]}\left|\int_r^T\left(p_n(s), q_n(s)\right)_H\,dW(s)\right|\leq
\delta\E\sup_{r\in[t,T]}\norm{p_n(t)}_H^{2} + C_\delta \E\norm{q_n}^{2}_{L^2(t,T; \cL^2(U,H))}\,,
\]
so that choosing $\delta$ sufficiently small and
using \eqref{est2_ad} we infer that 
\beq
  \label{est3_ad}
  \norm{p_n}_{L^{2}(\Omega; C^0([0,T]; H))}\leq c\,.
\eeq

\noindent{\bf Passage to the limit.}
We deduce that there is $(p,\tilde p,q)$ with 
\begin{gather*}
  p\in 
  L^\infty(0,T; L^{2}(\Omega; V))\cap L^{2}_{\mathcal P}(\Omega; L^2(0,T; Z\cap H^3(D)))\,,\\
  \tilde p \in 
  L^\infty(0,T; L^6(\Omega; V^*))\cap L^6_{\mathcal P}(\Omega; L^2(0,T; V))\,,\\
  q \in L^2_{\mathcal P}(\Omega; L^2(0,T; \cL^2(U,V)))\,,
\end{gather*}
such that $\tilde p=-\Delta p$ and
\begin{align*}
  p_n \wstarto p \qquad&\text{in } 
  L^\infty(0,T; L^{2}(\Omega; V))\cap L^{2}(\Omega; L^2(0,T; Z\cap H^3(D)))\,,\\
  \tilde p_n \wstarto \tilde p\qquad&\text{in } 
  L^\infty(0,T; L^6(\Omega; V^*))\cap L^6(\Omega; L^2(0,T; V))\,,\\
  q_n \wto q \qquad&\text{in } L^2(\Omega; L^2(0,T; \cL^2(U,V)))\,. 
\end{align*}
Now, we know from \cite[Lem.~2.1]{fland-gat} that the stochastic integral
operator is linear continuous (hence also weakly continuous)
from the space $L^2_{\mathcal P}(\Omega; L^2(0,T; \cL^2(U,V)))$ to 
the space
$L^2(\Omega; W^{s,2}(0,T; V))$: consequently, we deduce that 
\[
  \int_0^\cdot q_n(s)\,dW(s) \wto \int_0^\cdot q(s)\,dW(s) \qquad\text{in } 
  L^2(\Omega; W^{s,2}(0,T; V))\,.
\]
Finally, by (A2), the embedding $V\embed L^6(D)$ and the fact that
$y\in L^6(\Omega; L^\infty(0,T; V))$ 
it is immediate to check that $\Psi''(y) \in L^3(\Omega; L^\infty(0,T; L^3(D)))$, so 
in particular
\[
  \Psi''_n(y) \to \Psi''(y) \qquad\text{in } L^3(\Omega\times Q)\,.
\]
Hence, since by the convergences of $(\tilde p_n)_n$ we have
$\tilde p_n \wto \tilde p$ in $L^2(\Omega\times Q)$,
so that 
\[
  \Psi''_n(y)\tilde p_n \wto \Psi''(y)\tilde p \qquad\text{in } L^{6/5}(\Omega\times Q)\,.
\]
Similarly, it is a standard matter to check that 
the weak convergence of $(q_n)_n$ and the boundedness of $DB$ imply
\[
  DB(\cdot, y)^*q_n \wto DB(\cdot, y)^*q \qquad\text{in }
  L^2(\Omega; L^2(0,T; H))\,.
\]
Hence, passing to the weak limit as $n\to \infty$, we get that $(p,\tilde p, q)$
is a solution to the \eqref{ad1}--\eqref{ad4}.
Finally, note the extra regularities $p\in C^0_w([0,T]; L^{2}(\Omega; V))$ and 
$\tilde p \in C^0_w([0,T]; L^6(\Omega; V^*))$
follow a posteriori by comparison in the limit equation.

\subsection{Duality and conclusion}
In this final section we prove the last Theorem~\ref{thm:5} containing 
the simpler version of first-order necessary conditions on optimality.
The main idea is to remove the dependence on $z$ in the 
variational inequality 
of Theorem~\ref{thm:4} by using the adjoint problem and a 
suitable duality relation between $z$ and $\tilde p$.

Let then $\bar u \in \mathcal U$ be an optimal control and
$\bar y:=S(\bar u)$ be the corresponding solution to the state equation.
Then we know that the adjoint problem admits a solution $(p,\tilde p, q)$
solving \eqref{ad1}--\eqref{ad4}, where $\tilde p$ is uniquely determined.
Let $v\in \mathcal U$ be arbitrary
and set $h:=v-\bar u$:
the main point is to prove the duality relation
\[
  \alpha_1\E\int_Q(\bar y-x_Q)z_h + \alpha_2\E\int_D(\bar y(T) - x_T)z_h(T)
  =\E\int_Q\tilde p h \,.
\]
If we are able to prove such duality result, then it is clear that 
Theorem~\ref{thm:5} follows directly from Theorem~\ref{thm:4}.

Let $(z_h^n)_n$ and $(p_n, \tilde p_n, q_n)_n$ be the approximated solutions
introduced in Sections~\ref{subsec:lin} and \ref{subsec:ad}: then we have 
\begin{gather*}
  z_h^n \in L^2(\Omega; C^0([0,T]; H)\cap L^2(0,T; Z))\,,\\
  \tilde p_n \in C^0([0,T]; L^2(0,T; V^*))\cap L^2(\Omega; L^2(0,T; V))\,,
\end{gather*}
with $z_h^n(0)=0$, $p_n(T)=\alpha_2(\bar y(T)-x_T)$, and 
\begin{align*}
  &dz_h^n - \Delta (-\Delta z_h^n + \Psi''_n(\bar y)z_h^n - h)\,dt = 
  DB(\bar y)z_h^n\,dW\,,\\
  &-dp_n - \Delta\tilde p_n\,dt + \Psi''(\bar y)\tilde p_n\,dt = 
  \alpha_1(\bar y- x_Q)\,dt + D(\bar y)^*q_n\,dt - q_n\,dW\,,
\end{align*}
where the equations are intended in the Hilbert triplet $(Z,H,Z^*)$.
We deduce in particular that 
\[
  d(z_h^n, p_n)_H = (z_h^n,dp_n)_H + \ip{dz_h^n}{p_n}_Z +
  (DB(\bar y)z_h^n, q_n)_{\cL^2(U,H)}\,dt\,,
\]
where
\[
  (z_h^n,p_n)_H(0)=0\,, \qquad (z_h^n,p_n)_H(T)=\alpha_2\int_D(\bar y(T)-x_T)z_h^n(T)\,.
\]
Writing It\^o's formula for $(z_h^n, p_n)_H$ yields then
\begin{align*}
  &\alpha_2\E\int_D(\bar y(T)-x_T)z_h^n(T)
  =-\E\int_Q\Delta z_h^n \tilde p_n + \E\int_Q\Psi''_n(\bar y)\tilde p_n z_h^n
  -\alpha_1\E\int_Q(\bar y-x_Q)z_h^n\\
  &\quad-\E\int_0^T(DB(s, \bar y(s))^*q_n(s),z_h^n(s))_H\,ds
  -\E\int_Q\Delta z_h^n\Delta p_n
  +\E\int_Q\Psi''_n(\bar y) z_h^n\Delta p_n\\
  &\quad-\E\int_Q h\Delta p_n
  +\E\int_0^T\left(DB(s,\bar y(r))z_h^n(s), q_n(s)\right)_{\cL^2(U,H)}\,ds\,,
\end{align*}
from which, recalling the definition of $DB(\cdot, \bar y)^*$ and 
that $-\Delta p_n=\tilde p_n$, 
\[
  \alpha_1\E\int_Q(\bar y-x_Q)z_h^n +
  \alpha_2\int_D(\bar y(T)-x_T)z_h^n(T)=
  \E\int_Q h\tilde p_n \qquad\forall\,n\in\enne\,.
\]
The thesis now follows letting $n\to\infty$.

\bibliographystyle{abbrv}
\def\cprime{$'$}

\end{document}